\documentclass[12pt]{amsart}
\usepackage[active]{srcltx}
\begin{document}
\numberwithin{equation}{section}

\def\1#1{\overline{#1}}
\def\2#1{\widetilde{#1}}
\def\3#1{\widehat{#1}}
\def\4#1{\mathbb{#1}}
\def\5#1{\frak{#1}}
\def\6#1{{\mathcal{#1}}}

\def\C{{\4C}}
\def\R{{\4R}}
\def\N{{\4N}}
\def\Z{{\4Z}}

{\em American Journal of Mathematics}, {\bf 134} (2012), no. 4, 915--947.
\bigskip

Published by Johns Hopkins University Press

\medskip

DOI: https://doi.org/10.1353/ajm.2012.0027

https://muse.jhu.edu/article/481304/summary

\bigskip
\bigskip
\bigskip
\bigskip

\title[Normal forms for nonintegrable almost CR structures]{Normal forms for nonintegrable almost CR structures}
\author[D. Zaitsev]{Dmitri Zaitsev*}
\thanks{*Supported in part by the Science Foundation Ireland grant 06/RFP/MAT018.}
\address{D. Zaitsev: School of Mathematics, Trinity College Dublin, Dublin 2, Ireland}
\email{zaitsev@maths.tcd.ie}

\begin{abstract}
We propose two constructions extending the Chern-Moser normal form to non-integrable Levi-nondegenerate (hypersurface type) almost CR structures. One of them translates the Chern-Moser normalization into pure intrinsic setting, whereas the other directly extends the (extrinsic) Chern-Moser normal form by allowing non-CR embeddings that are in some sense ``maximally CR". One of the main differences with the classical integrable case is the presence of the non-integrability tensor at the same order as the Levi form, making impossible a good quadric approximation - a key tool in the Chern-Moser theory. Partial normal forms are obtained for general almost CR structures of any CR codimension, in particular, for almost-complex structures. Applications are given to the equivalence problem and the Lie group structure of the group of all CR-diffeomorphisms. 
\end{abstract}

\maketitle
\tableofcontents

\def\Label#1{\label{#1}}


\def\cn{{\C^n}}
\def\cnn{{\C^{n'}}}
\def\ocn{\2{\C^n}}
\def\ocnn{\2{\C^{n'}}}


\def\dist{{\rm dist}}
\def\const{{\rm const}}
\def\rk{{\rm rank\,}}
\def\id{{\sf id}}
\def\tr{{\bf tr\,}}
\def\aut{{\sf aut}}
\def\Aut{{\sf Aut}}
\def\CR{{\rm CR}}
\def\GL{{\sf GL}}
\def\Re{{\sf Re}\,}
\def\Im{{\sf Im}\,}
\def\span{\text{\rm span}}
\def\Diff{{\sf Diff}}

\def\codim{{\rm codim}}
\def\crd{\dim_{{\rm CR}}}
\def\crc{{\rm codim_{CR}}}

\def\phi{\varphi}
\def\eps{\varepsilon}
\def\d{\partial}
\def\a{\alpha}
\def\b{\beta}
\def\g{\gamma}
\def\G{\Gamma}
\def\D{\Delta}
\def\Om{\Omega}
\def\k{\kappa}
\def\l{\lambda}
\def\L{\Lambda}
\def\z{{\bar z}}
\def\w{{\bar w}}
\def\Z{{\1Z}}
\def\t{\tau}
\def\th{\theta}

\emergencystretch15pt
\frenchspacing

\newtheorem{Thm}{Theorem}[section]
\newtheorem{Cor}[Thm]{Corollary}
\newtheorem{Pro}[Thm]{Proposition}
\newtheorem{Lem}[Thm]{Lemma}

\theoremstyle{definition}\newtheorem{Def}[Thm]{Definition}

\theoremstyle{remark}
\newtheorem{Rem}[Thm]{Remark}
\newtheorem{Exa}[Thm]{Example}
\newtheorem{Exs}[Thm]{Examples}

\def\bl{\begin{Lem}}
\def\el{\end{Lem}}
\def\bp{\begin{Pro}}
\def\ep{\end{Pro}}
\def\bt{\begin{Thm}}
\def\et{\end{Thm}}
\def\bc{\begin{Cor}}
\def\ec{\end{Cor}}
\def\bd{\begin{Def}}
\def\ed{\end{Def}}
\def\br{\begin{Rem}}
\def\er{\end{Rem}}
\def\be{\begin{Exa}}
\def\ee{\end{Exa}}
\def\bpf{\begin{proof}}
\def\epf{\end{proof}}
\def\ben{\begin{enumerate}}
\def\een{\end{enumerate}}
\def\beq{\begin{equation}}
\def\eeq{\end{equation}}

\section{Introduction}\Label{intro}
The main goal of this paper is to extend the Chern-Moser normal form theory \cite{CM}
to {\em non-integrable} (or {\em non-involutive}) almost CR structures.
Recall that an {\em almost CR structure} on a real manifold $M$
consists of a subbundle $H=HM$ of the tangent bundle 
$T=TM$ and a vector bundle automorphism $J\colon H\to H$ satisfying $J^2=-\id$.
Thus $J$ makes every fiber $H_p$, $p\in M$, into a complex vector space.
Two important special cases of 
almost CR structures are {\em almost complex structures} corresponding to $H=T$ and {\em hypersurface type almost CR structures} corresponding to $H$ being corank $1$ subbundle of $T$. Equivalently, an almost CR structure is given by the $i$-eigenspace subbundle $H^{1,0}\subset\C\otimes H\subset\C\otimes T$ of (the  $\C$-linear extension of) $J$, which can be chosen as an arbitrary complex subbundle of $\C\otimes T$ satisfying $H^{1,0}\cap \1{H^{1,0}}=0$. A {\em CR structure} is an almost CR structure satisfying the integrability condition $[H^{1,0},H^{1,0}]\subset H^{1,0}$. The reader is referred to \S\ref{notions} for further basic notions for almost CR structures used throughout the paper.

The classification problem for CR structures goes back to H.~Poincar\'e~\cite{Po}
who realized that the latters have infinitely many local invariants and posed the problem
of finding those invariants. Solutions to this problem for Levi-nondegenerate
hypersurface type CR structures were given by
E.~Cartan~\cite{CaE} in real dimension $3$ and by  N. Tanaka~\cite{Ta62} and S.S. Chern and J. Moser \cite{CM} in the general case. One of the main tools provided in \cite{CM}
 for classification of CR structures and finding their invariants
 is a {\em normal form}. The latter is a realization of every real-analytic Levi-nondegenerate hypersurface type CR structure near a given point $p$ as an embedded real-analytic hypersurface $M\subset\C^{n+1}$ of the form 
\begin{equation}\Label{graph}
\Im w= F(z,\bar z,\Re w),  \quad (z,w)\in\C^n\times\C,
\end{equation}
 with $p$ corresponding to the origin and $F$ having the expansion
 \begin{equation}\Label{cm-expansion}
 F(z,\bar z,u)=\sum_{j=1}^n \eps_j |z_j|^2+ 
\sum_{k,l\ge2} F_{kl}(z,\bar z,u),
\end{equation}
where 
\begin{equation}\Label{eps}
\eps_j=
\begin{cases}
-1 & j\le s\cr
1 & j>s
\end{cases}
\end{equation}
 for suitable $s$
(the signature of the Levi form)
and the summands $F_{kl}(z,\bar z,u)$ are bihomogeneous of bidegree $(k,l)$ in $(z,\bar z)$
(i.e.\ $F_{kl}(tz,r\bar z,u)=t^k r^l F_{kl}(z,\bar z,u)$ for $t,r\in\R$)
 satisfying the normalization conditions
\begin{equation}\Label{cm-normalization}
\tr F_{22}\equiv0, \quad \tr^2 F_{23}\equiv0, \quad \tr^3 F_{33}\equiv0,
\end{equation}
where $\tr : = \sum_{j=1}^n \eps_j \frac{\partial^2}{\partial z_j \partial \bar z_j}$.
The CR structure here is induced by the embedding of $M$ into $\C^{n+1}$ in the sense that $HM=TM\cap iTM$ and $J$ is the restriction to $HM$ of the multiplication by $i$ in $\C^{n+1}$.
If the given CR structure is merely smooth, it still admits a {\em formal normal form} \eqref{graph} at every point with $F$ being a formal power series satisfying \eqref{cm-expansion} and \eqref{cm-normalization}. 
(This is a consequence of the formal CR embeddability, cf.\ e.g.\ \cite[Proposition~3.1]{KZ} and \cite[Theorem~2.1.11]{BERbook}, and the formal normal form \cite[Theorem~2.2]{CM}.)
The latter fact explains  the significance of
{\em formal (almost) CR structures}, i.e.\
those given by formal power series,
that will be considered throughout this paper.  
More precisely, a {\em formal almost CR structure}
on a smooth manifold $M$ at a point $p\in M$
is given in suitable smooth local coordinates $(x,y)\in\R^a\times\R^b$
on $M$ vanishing at $p$ by a {\em formal subbundle}
$H=\{dy=h(x,y)dx\}$ and a {\em formal endomorphism}
$J=j(x,y)\colon\R^a\to\R^a$,
where $h(x,y)$ and $j(x,y)$ are matrix-valued formal power series
satisfying $j(x,y)^2=-\id$.

 Our goal here is to obtain a normal form of the above type for more general almost CR structures that may not necessarily satisfy the integrability condition.
 Such structures arise naturally, for instance, when one is deforming or glueing CR structures. Another source of almost CR structures is given by real submanifolds of general almost complex manifolds.
 
 The first problem one faces when extending the Chern-Moser normal form approach
is that the non-integrable almost CR structures do not admit any realization as real submanifolds in $\C^{n+1}$ even at the formal level. 
Speaking informally, the normal form 
\eqref{graph}--\eqref{cm-normalization} is {\em extrinsic}
(i.e.\ written in the ambient space coordinates of $\C^{n}\times\C$) whereas almost CR structures are {\em intrinsic}
(i.e.\ written in local coordinates of the manifold $M$ itself).

In this paper we suggest two different ways of overcoming this difficulty.
The first way is based on an ``intrinsic analogue'' of the Chern-Moser normal form.
The starting idea is to impose normalization conditions on the almost CR structure itself rather than on the defining equation of an embedding.
However, almost CR structures are given by objects of different nature (complex subbundles of $\C\otimes T$) than defining functions.
In order to relate both settings we introduce a new function $F$ associated
with a given almost CR structure in given local coordinates on the manifold. Roughly speaking, $F$ is obtained by ``restricting" the almost CR structure to the Euler vector field. 
More precisely, consider local coordinates on $M$
that we group as $(z,u)\in\C^n\times\R$, where the subspace
$\C^n\times\{0\}$ with its complex structure
corresponds to the given almost CR structure on $M$ at the origin (but not necessarily at other points).
Using the standard complex structure we write the complexification of $\C^n$
as the direct sum $\C^n_z\oplus\C^n_{\bar z}$ of the spaces of $(1,0)$ and $(0,1)$ vectors
(i.e.\ $\pm i$-eigenspaces of $J$). Then $\C\otimes T$ can be identified  at every point
with $\C^n_z\oplus\C^n_{\bar z}\oplus\C_w$ with $u=\Re w$ and the almost CR structure is given by a complex subbundle $H^{1,0}\subset \C\otimes T$
which, at every point $p=(z,u)$ near $0$, is the graph of a uniquely determined complex-linear
map $L(z,\bar z,u)\colon \C^n_z\to\C^n_{\bar z}\oplus\C_w$.
We consider the Euler (or radial) vector field $e(z)=\sum_j z_j\frac{\d}{\d z_j}$ on $\C^n_z$ and define the ``radial part" of $L$ as the evaluation along $e(z)$, i.e.\
\begin{equation}\Label{L}
\2L(z,\bar z,u):=L(z,\bar z,u)(e(z)), \quad \2L=(\2L^{\bar z},\2L^w)\in \C^n_{\bar z}\oplus\C_w.
\end{equation}

Now, in order to write normalization conditions similar to the ones by Chern-Moser,
we consider the Taylor series expansions  in $(z,\bar z,u)$ at $0$ of the functions $L$ and $\2L$.
Since we consider arbitrary smooth almost CR structures,
the corresponding Taylor series are {\em formal}.
By a slight abuse of notation, we denote these Taylor series
by the same letters $L$ and $\2L$ respectively.
Our goal now is to obtain a {\em formal normal form} for $L$,
i.e.\ a set of normalization conditions for the coefficients of $L$
achieved by a suitable {\em formal change of coordinates} on $M$.
We first state our main result for almost CR structures with 
{\em positive definite
Levi form} and later explain how it can be extended to more general nondegenerate Levi forms.

\bt\Label{main-intrinsic}
Any (germ of a) smooth hypersurface type almost CR structure with positive definite Levi form admits a formal normal form in $\C^n_z\times\R_u$ given by a formal map 
$L(z,\bar z,u)\colon \C^n_z\to\C^n_{\bar z}\oplus\C_w$ with
\begin{equation}\Label{norm-intrinsic}
\2L^{\bar z}\equiv0,  \quad \Re \2L^w\equiv0,
\end{equation}
 and $F:=\Im\2L^w$ satisfying the Chern-Moser normalization conditions
\eqref{cm-expansion}--\eqref{cm-normalization}.
\et

Theorem~\ref{main-intrinsic} is a direct consequence of Theorem~\ref{full} below.
An important feature of the normal form in Theorem~\ref{main-intrinsic}
is that different normal forms of the same almost CR structure
are determined precisely by the same set of parameters as the 
Chern-Moser normal form, i.e.\ by the automorphism group 
of the associated hyperquadric.
 cise statement is given in Theorem~\ref{uniqueness} below.

Conditions~\eqref{norm-intrinsic} admit a natural geometric interpretation.
The tangential component $L^{\bar z}$ always automatically vanishes identically
for the CR structure of a hypersurface of the form \eqref{graph}, where we consider
$z$ and $u=\Re w$ as intrinsic coordinates on the hypersurface.
Now the first condition in \eqref{norm-intrinsic}
expresses the property that
the ``radial part" $\2L^{\bar z}$ of the tangential component $L^{\bar z}$
can still be eliminated in general.
Again, going back to $\eqref{graph}$, we see that the function $F$ corresponds to the imaginary part of $w$ and is therefore automatically real. 
On the other hand, the ``radial part" $\2L^w$ of the transversal component $L^{w}$ is complex-valued in general.
Then the second condition in \eqref{norm-intrinsic}
asserts that its real part can be completely eliminated
and thus $\2L^w$ is given by its imaginary part $F$.
The presence of conditions~\eqref{norm-intrinsic} 
reflects the fact that the moduli space of almost CR structures
has additional parameters compared to the moduli space of (integrable) CR structures,
where the purpose of conditions~\eqref{norm-intrinsic}
is precisely to restrict those additional parameters.

As a byproduct we also obtain a {\em partial normal form} for almost CR structures
that are of {\em arbitrary CR codimension} (not necessarily of hypersurface type)
and {\em without any nondegeneracy conditions}, see Proposition~\ref{partial} below.
It can be seen as an extension to the non-integrable case of the so-called {\em normal coordinates}
(see e.g.\ \cite{BERbook}). Proposition~\ref{partial} also describes the degree of uniqueness
of the partial normal form. Similarly to the normal form in Theorem~\ref{main-intrinsic}
being parametrized by the CR-automorphisms of the hyperquadric, the partial normal form
in Proposition~\ref{partial} is parametrized by the CR-automorophisms of the standard CR structure 
on $\C^n\times\R^d$ with $H=T\C^n\times\{0\}$, 
i.e.\ by (formal) maps $(z,u)\mapsto (f(z,u),g(u))$ preserving the origin
with $f(z,u)$ being holomorphic in $z$.

In particular, if the CR codimension is zero, i.e.\ if we are given
an {\em almost-complex structure}, it admits at every point a (formal) normal form
with vanishing $\2L^{\bar z}$. The latter condition is equivalent to the almost CR structure
being the standard one along the complex lines through the origin.
That is, we obtain a normal form for arbitrary almost-complex structures.
This normal form is
characterized by the property that
all complex lines through the origin are pseudo-holomorphic
(i.e.\ the map $t\in\C\mapsto tv\in\C^n$  is holomorphic for every $v\in\C^n$ with the respect
to the almost-complex structure on $\C^n$ in its normal form.
Specializing the above uniqueness remark to the case of CR codimension zero,
we conclude that the normal form in this case is parametrized by (formal)
biholomorphic self-maps of $\C^n$ preserving the origin.
(See also \cite{Kr1, Kr2,To} for other normal forms for almost-complex structures.)

The normal form given by Theorem~\ref{main-intrinsic} is intrinsic and mimics the Chern-Moser normalization.
However, specializing to CR structures (i.e.\ the integrable ones), 
one obtains here an intrinsic normal form that is different from the (extrinsic) Chern-Moser normal form of the same CR structure.
The goal of our second construction is to obtain a normal form that directly extends the Chern-Moser one. The idea is to go back to the extrinsic point of view, using the ambient space coordinates in $\C^{n+1}$, where $M$ is embedded, but to allow embeddings that are not CR. That is, instead of looking for CR embeddings that may not exist, we consider certain more general embeddings called here {\em quasi CR embeddings} that always exist and are in certain sense ``maximally CR". More precisely, a {\em quasi CR embedding} of $M$ at a point $p$ in $\C^{n}_z\times\C_w$ is an embedding as a hypersurface $M'$ with $p=0$, $T_0M'=\C^n\times\R$, and such that the transformed almost CR structure on $M'$ 
and the one induced by the embedding have equal ``radial parts''
in the sense of \eqref{L}, i.e.\ $\2L\equiv\2L'$
for the corresponding (Taylor series expansions of the) maps $L$ and $L'$,
where we regard $z$ and $u=\Re w$ as coordinates on $M'$.
The latter property can be rephrased by saying that both almost CR
structures coincide along the Euler vector field $e(z)$.

\bt\Label{main-quasi}
Any (germ of a) smooth hypersurface type almost CR structure
 with positive definite Levi form admits a normal form consisting of a formal quasi CR embedding into $\C^{n+1}$
as a formal real hypersurface satisfying the Chern-Moser normalization \eqref{cm-expansion}--\eqref{cm-normalization}.
A (formal) quasi CR embedding is an actual CR embedding if and only if the given almost CR structure is integrable, in which case this normal form coincides with the Chern-Moser one.
\et

As with the normal form in Theorem~\ref{main-intrinsic},
the one given by Theorem~\ref{main-quasi} has the property that
different normal forms of the same almost CR structure
depend on the same parameters as the Chern-Moser normal form
(see Theorem~\ref{uniqueness} below). 
Also here we obtain a partial normal form for almost CR structures
of any codimension without any nondegeneracy assumptions, see Proposition~\ref{quasi-prop}.
It states that any almost CR structure admits a quasi CR embedding 
as a generic submanifold in normal coordinates (see \cite{BERbook}).
This is another extension of normal coordinates to the non-integrable case
with Remark~\ref{forms-different} below showing that
the two extensions given by Propositions~\ref{partial} and \ref{quasi-prop}
are actually different. On the other hand, in the case of almost-complex structures
(CR codimension zero) both normal forms are equal.

Let us now mention the crucial difference between the situations of Theorems~\ref{main-intrinsic} - \ref{main-quasi} and the integrable case of Chern-Moser \cite{CM}. In the latter, the Levi form is the only second order obstruction for a CR structure from being {\em flat}, i.e.\ induced by a hyperplane embedding in $\C^m$. Consequently, any CR structure can be approximated to the third order
by the hyperquadric associated with its Levi form (even to the forth order if the Levi form is nondegenerate). Chern-Moser approach is heavily based on this approximation.
However, in the non-integrable case, it is no more true that the Levi form is the only lowest order obstruction to the flatness (i.e.\ being induced by a hyperplane embedding). Indeed, the non-integrability (Nijenhuis) tensor appears at the same order and ``interferes" with the Levi form in the sense that it arises everywhere in the transformation formulas along with the Levi form terms. 
It is no more possible to separate the ``good" Levi form terms from the rest by looking only at their weights. As a consequence, the decoupled transformation formulas in Chern-Moser case, become non-decoupled in our case and the order in which the normalization 
proceeds, becomes essential.
The presence of the non-integrability tensor also prevents almost CR structures from being approximable to the third order by hyperquadrics. 
It is therefore remarkable that one can still extend the Chern-Moser approach as in Theorems~\ref{main-intrinsic} - \ref{main-quasi} with assumptions put only on the Levi form.

The mentioned difficulties become even more apparent when one tries to extend Theorems~\ref{main-intrinsic}-\ref{main-quasi} to the case of  the Levi form being {\em nondegenerate}
rather than positive (or negative) definite. Here, in addition to the nondegeneracy of the Levi form itself, 
one also needs the nondegeneracy of
a certain linear combination of the Levi form and the transversal component of the non-integrability tensor. More precisely, given the function 
$L = (L^{\bar z}, L^w)\colon \C^n_z\to\C^n_{\bar z}\oplus\C_w$
as above defining the almost CR structure, the Levi form at $0$ corresponds, up to an imaginary multiple, to the antihermitian part $\6L(\bar\xi,\eta)$ of the derivative of $L^w$
in ${\bar z}$, whereas the non-integrability tensor corresponds to the antisymmetric part
$\6N(\xi,\eta)$ of the derivative of $L$ in $z$, and its transversal component $\6N^w(\xi,\eta)$
to that of $L^w$. 
(In fact, $\6L$ and $\6N$ are the only 2nd order obstructions to the flatness
as immediately follows from our normal form.)
We now call the almost CR structure {\em strongly nondegenerate} at $0$
if it is Levi-nondegenerate and in addition, the bilinear form $6i\6L+\6N^w$ is nondegenerate in the sense that
\begin{equation}\Label{str-nondeg}
6i\6L(\bar\xi,\eta)+\6N^w(\xi,\eta)=0 \text{ for all } \eta \quad  \Longrightarrow
\quad \xi=0.
\end{equation}
Obviously, for (integrable) CR structures (i.e.\ with $\6N=0$), strong nondegeneracy means the same as Levi-nondegeneracy. Furthermore, a strongly pseudoconvex almost CR structure is automatically strongly nondegenerate. Indeed, since $\6N$ is antisymmetric, substituting $\eta=\xi$ into the left-hand side of \eqref{str-nondeg} leads to $\6L(\bar\xi,\xi)=0$, which in the case
$\Im\6L$ is positive definite, implies $\xi=0$. However, if the Levi form has mixed signature, strong nondegeneracy is a stronger property than Levi-nondegeneracy. 
We now have the following extensions of Theorems~\ref{main-intrinsic} and \ref{main-quasi}
covering, in particular, all Levi-nondegenerate (integrable) CR structures:

\bt\Label{strong-both}
Conclusions of Theorems~\ref{main-intrinsic} and \ref{main-quasi}
hold for strongly nondegenerate hypersurface type almost CR structures.
\et

Theorem~\ref{strong-both} is a direct consequence of Theorems~\ref{full} and \ref{full-extrinsic} below.
We conclude by describing the uniqueness of the above normal forms. 
As mentioned before, each of the normal forms of the same almost CR
structure is determined by the same parameters as the one of Chern-Moser.
The latter is known to be parameterized by the automorphism group of the associated hyperquadric. However, this parametrization is not unique.
Instead, we take here a more direct geometric approach
describing naturally the needed parameters.

First recall that an {\em adapted frame} of an
(almost) CR structure on $M$ at $p\in M$
with non-degenerate Levi form $\6L_p\colon H_p\times H_p\to \C\otimes(T_p/H_p)$
of signature $s$ consists of a complex basis $v_1,\ldots,v_n$ of $H_pM$,
extended by a vector $v_{n+1}\in T_pM\setminus H_pM$,
satisfying
\begin{equation}\Label{orthog}
\6L_p(v_j,v_k)= \eps_j\delta_{jk}[v_{n+1}],
\end{equation}
where $\eps_j$ is as in \eqref{eps},
$\delta_{jk}$ is Kronecker delta
and $[v_{n+1}]\in T_p/H_p$ stands for the equivalence class of $v_{n+1}$.
In given coordinates $(x,y,u)\in\R^n\times\R^{n}\times\R$ on $M$ with $H_0M=\C^n\times\{0\}$,
one has the {\em standard adapted frame}
$\frac{\d}{\d x_1},\ldots, \frac{\d}{\d x_n}, \frac{\d}{\d u}$ at $0$.
It follows from Chern-Moser theory \cite{CM} that 
any Levi-nondegenerate CR structure can be always be put into 
its Chern-Moser normal form with the additional
condition that an arbitrary given adapted frame is transformed into
the standard one.
However, the latter condition still does not suffice to make the normal form unique.
In order to achieve the uniqueness, one needs a second order condition,
of which we give here the following geometric version.
We here give a brief description and refer the reader to \S\ref{jets} for more details.

Given and adapted frame $v_1,\ldots,v_n,v_{n+1}$ at $p$,
consider $2$-jets $\L$ at $0$ of real curves in $M$, 
given by equivalence classes of smooth curve germs $\gamma\colon(\R,0)\to (M,p)$
with $\gamma(0)=p$ and $\gamma'(0)=v_{n+1}$.
We say that two such jets $\L_1$ and $\L_2$ are {\em $H$-equivalent}
if 
\begin{equation}\Label{h-eq}
\gamma''_1(0)-\gamma''_2(0)\in H_pM
\end{equation}
holds
for some represenatives $\gamma_1$, $\gamma_2$ of $\L_1$, $\L_2$
respectively and where the second derivatives are calculated in some local
coordinates on $M$. It follows from the canonical affine structure
on jet bundles that the property \eqref{h-eq}
is independent of both representatives as well as of the coordinate choice involved
(see \S\ref{jets}). Thus $H$-equivalence is a well-defined equivalence relation.
We now define an {\em extended adapted frame} on $M$ at $p$
to be any collection $v_1,\ldots,v_n,v_{n+1},[\L]$,
consisting of an adapted frame together with
a choice of an $H$-equivalence class of $2$-jets $\L$ as above.
Given coordinates $(x,y,u)$ as above,
one can naturally complete the standard adapted frame
$\frac{\d}{\d x_1},\ldots, \frac{\d}{\d x_n}, \frac{\d}{\d u}$
to the {\em standard extended adapted frame}
by taking $[\L]$ with $\L$ represented by 
the curve $\gamma(t)=(0,t)$.
It now follows from \cite{CM} that any Levi-nondegenerate
(hypersurface type)
CR structure can always be put into its Chern-Moser normal form
with the additional condition that
an arbitrary given extended adapted frame
is transformed into the standard one.
This time the normal form is uniquely determined.
Thus the space of all extended adapted frames
parametrizes precisely the space of all Chern-Moser
normal forms of a given CR structure.
The main uniqueness result for our normal forms 
states that the latters are parametrized in the same way:

\bt\Label{uniqueness}
Let $(M,p)$ be a germ of a smooth real manifold with strongly nondegenerate hypersurface type almost CR structure.
Then for every extended adapted frame at $p$,
each of the normal forms in Theorems~\ref{main-intrinsic} and \ref{main-quasi}
exists and is unique under the additional condition
that the given extended adapted frame is realized as the standard one in the normal form.
\et

Theorem~\ref{uniqueness} is a direct consequence of Theorem~\ref{unique-precise} below.

\section{Applications of the normal forms}

We here state some application of our normal form.
The statements are formulated for strongly nondegenerate CR structures.
Recall that any hypersurface type almost CR structure with {\em positive definite Levi form} is 
automatically strongly nondegenerate.
Also any hypersurface type (integrable) CR structure with {\em nondegenerate Levi form} is
automatically strongly nondegenerate. 

Given two manifolds $M$ and $M'$ of the same dimension and an integer $k$, 
we denote by $G^k(M,M')$ the space of all $k$-jets of invertible maps. 
For a point $p\in M$, denote by $G^k_p(M,M')\subset G^k(M,M')$
the subset of jets with source $p$.
For any smooth diffeomorphism $f$ between open pieces of $M$ and $M'$
and any $p\in M$ in the domain of definition of $f$, we denote by $j^k_p f\in G^k(M,M')$
the $k$-jet of $f$ at $p$.

Recall that a {\em CR-diffeomorphism} between two almost CR manifolds
is any diffeomorphism transforming the almost CR structure of the first manifold
into the structure of the second.
Our first application is given by the following complete system type property of local CR-diffeomorphism.

\bt\Label{system}
Let $M$ and $M'$ be smooth (resp.\ real-analytic) strongly nondegenerate hypersurface type almost CR manifolds of the same dimension. 
Then for every $p_0\in M$ and $\L_0\in G^2_{p_0}(M,M')$, there exist
an open neighborhood $\Omega$ of $\L_0$ in $G^2(M,M')$ and  a 
smooth (resp.\ real-analytic) map $\Phi\colon \Omega\to G^3(M,M')$ such that
every smooth CR-diffeomorphism $f$ between open pieces of $M$ and $M'$
satisfies the complete differential system
\begin{equation}\Label{sys}
j^3_pf = \Phi(j^2_pf)
\end{equation}
whenever $p\in M$ is in the domain of definition of $f$ and $j^2_pf\in \Omega$.
\et

As an application of Theorem~\ref{system}
we obtain the following unique determination result.

\bc\Label{unique}
Let $M$ and $M'$ be as in Theorem~\ref{system}.
Then CR-diffeomorphisms between connected open pieces of $M$ and $M'$
are uniquely determined by their $2$-jets at any point of $M$,
i.e. if $f,g\colon U\subset M\to M'$ are CR-diffeomorphisms, $U$ is connected,
$p\in U$ and $j^2_p f=j^2_p g$, then $f\equiv g$.
\ec

In fact, we have the following stronger property
that local CR-diffeomorphisms depend
smoothly (resp.\ analytically) on their $2$-jets.

\bc\Label{param}
Let $M$ and $M'$ be as in Theorem~\ref{system}.
Then for every $p_0\in M$ and
$\L_0\in G^2_{p_0}(M,M')$, there exist   neighborhoods
$\Omega'$ of
$p_0$ in $M$,
$\Omega''$ of $\L_0$  in  $G^2(M,M')$,
  and a smooth
(resp. real-analytic) map
$\Psi\colon\Omega'\times \Omega''\to M'$ 
such that 
every smooth CR-diffeomorphism $f$ between open pieces of $M$ and $M'$
satisfies
\begin{equation}\Label{jps}
f(q)
=
\Psi(q,j^{2}_p f),
\end{equation}
whenever  $q\in \Omega'$ and $p\in M$ are in the domain of definition of $f$ and
 $j^{2}_p f\in
\Omega''$.
In particular, if $M$ and $M'$ are real-analytic, any smooth CR-diffeomorphism 
between them is also real-analytic.
\ec

The proofs of Corollaries~\ref{unique} and \ref{param} are completely analogous to that of Proposition~2.2 in \cite{BRWZ} (see also \cite{Ejets} and \cite{KZ}).
Corollary~\ref{param} can be used to obtain a solution to the CR equivalence
problem in the following sense:

\bt\Label{loc-equiv}
Let $M$ and $M'$ be as in Theorem~\ref{system}.
Then for every $p_0\in M$ and
$\L_0\in G^2_{p_0}(M,M')$, there exist an open neighborhood
$\Omega$ of $(p_0,\L_0)$ in $M\times G^2(M,M')$ 
and a smooth (resp.\ real-analytic) real function $\Gamma\colon\Omega\to \R$
such that there exists a local CR-diffeomorphism $f$
between an open neighborhood of a point $p\in M$ and an open set in $M'$
with $(p,j^2_pf)=(p,\L)\in \Omega$
if and only if  $\Gamma(\cdot,\L)$ vanishes
in a neighborhood of $p$.
\et

The proof is completely analogous to that of Theorem~1.4 in \cite{KZ}.
Finally, as a direct consequence of Theorem~\ref{system} above and Theorem~2.3 in \cite{BRWZ}, we obtain a Lie group structure 
on the group of all CR-automorphisms in its natural topology.
Recall that the group $\Diff(M)$ of all smooth diffeomorphisms of $M$
carries the natural topology of uniform convergence on compacta of maps
and their inverses together with all derivatives of the maps and the inverses.
We call this topology {\em  the compact-open $C^\infty$ topology}.
Applying Theorem~\ref{sys} and the mentioned result from \cite{BRWZ}, we obtain:

\bc\Label{main-lie} 
Let $M$ be a smooth connected
strongly nondegenerate hypersurface type almost CR manifold.
Then the group of all smooth CR-automorphisms of $M$, 
equipped with the compact-open $C^\infty$ topology, 
has a (unique) structure as a Lie group
 acting smoothly on $M$.
\ec

In case $M$ is real-analytic, the group of all real-analytic CR-automorphisms of $M$
(which coincides with the group of all smooth CR-automorphisms in view of Corollary~\ref{param}) can also be equipped with the natural compact-open real-analytic topology
(see \cite{BRWZ} for details).
Then Theorem~\ref{sys} can be applied along with Theorem~2.3 in \cite{BRWZ}
to obtain the real-analytic analogue of Corollary~\ref{main-lie}.

\section{Almost CR structures}\Label{notions}
Let $M$ be a real manifold with almost CR structure given 
by $H$ and $J$ or, equivalently by a complex subbundle $H^{1,0}\subset \C\otimes T$ 
satisfying $H^{1,0}\cap \1{H^{1,0}}=0$.
The complex dimension of the fiber $H_p$, $p\in M$, is called the {\em CR dimension} 
$\crd M$ of
an almost CR manifold $M$ and the real codimension of 
$H_p$ in $T_p$ the {\em CR codimension} $\crc M$ of $M$.
The pair $(\crd M,\crc M)$ is sometimes called the {\em type} of $M$
and in case $\crc M=1$, $M$ is said to be of {\em hypersurface type}.

Many formulas and calculations become simpler when working with
the complexified tangent bundle $\C T:= \C\otimes T$.
(For instance, the Nijenhuis tensor of the almost complex structure
is expressed by a linear combination of $4$ different brackets of real vector fields but only one bracket of complexified vector fields.)
Here $J$ extends to a complex bundle automorphism of 
$\C H=:\C\otimes H$, which splits into direct sum
of its $(\pm i)$-eigenspaces 
$$H^{1,0}:= \{\xi\in \C H : J\xi = i\xi\}, 
\quad H^{0,1}:= \{\xi\in \C H : J\xi = -i\xi\}.$$
These eigenspaces form complex subbundles of $\C H$ satisfying
$$H^{0,1}=\1{H^{1,0}},\quad \C H = H^{1,0} \oplus H^{0,1},$$
and each of them uniquely determines the almost CR structure.

\subsection{The non-integrability tensor and the Levi form}
Recall that an {\em almost CR structure} is called 
{\em (formally) integrable}
or simply a {\em CR structure} if the subbundle $H^{1,0}\subset \C T$
(or equivalently $H^{0,1}$) is closed under taking Lie brackets.
Denote by 
$\pi^{0,1}_p\colon \C T_p \to \C T_p / H^{1,0}_p$
the canonical projection.
Then the obstruction to integrability is given by the antisymmetric
complex-bilinear map
\begin{equation}\Label{N}
\6N_p\colon H^{1,0}_p \times H^{1,0}_p \to \C T_p / H^{1,0}_p
\end{equation}
such that 
$\6N_p(X_p,Y_p)= \pi^{0,1}_p([X,Y]_p)$
holds for any vector fields $X,Y\in \Gamma(H^{1,0})$,
where $X_p$ denotes the evaluation at $p$.
Then the CR structure is formally integrable if and only if $\6N_p=0$ 
for all $p\in M$.
The map $\6N_p$ generalizes the {\em Nijenhuis tensor}
of almost complex structures to arbitrary almost CR structures
and is called here the {\em non-integrability tensor}. 

The {\em Levi form} is defined for any almost CR structure in a similar fashion.
This time consider the canonical projection
 $\pi_p\colon \C T_p\to \C T_p/\C H_p$.
Then the Levi form is given by the sesqui-linear map
$$\6L^{1,0}_p\colon H^{1,0}_p\times H^{1,0}_p \to \C T_p/\C H_p$$
such that $\6L^{1,0}_p(X_p,Y_p)=\frac1{2i}\pi_p([\bar X,Y]_p)$
holds for any vector fields $X,Y\in \Gamma(H^{1,0})$
(cf.\ e.g.\ \cite{BERbook}).
Using the natural identifications $H^{0,1} \to H^{1,0}$ given by 
$\xi\mapsto\bar\xi$, and $H\to H^{1,0}$ given by
$\xi\mapsto \frac12(\xi-iJ\xi)$,
one obtains the corresponding maps 
$$\6L^{0,1}_p\colon H^{0,1}_p\times H^{0,1}_p \to \C T_p/\C H_p,
\quad \6L_p\colon H_p\times H_p \to \C T_p/\C H_p,$$
of which the first is used in \cite{BERbook}
and the second one will be used here.
(With this normalization the Levi form of the quadric 
$\Im w=|z|^2$ gives $\6L_0(s\d_x,t\d_x)=\bar s t\d_u$,
where $z=x+iy$, $u=\Re w$, and $\C T_p/\C H_p$ is identified with $\C_w$.)

\section{Coordinate setting}
\subsection{Preliminary normalization}
Given an almost CR structure $(H,J)$ on $M$ of CR dimension $n$ and
CR codimension $d$ and a reference point $p_0\in M$,
we consider a system of coordinates
$(z,u)=(x+iy,u)\in\C^n\times\R^d$
on $M$ such that
\begin{equation}\Label{preli-norm}
p_0=(0,0,0), \quad H_0=\C^n\times\{0\}, \quad J_0(\xi,0)=(i\xi,0).
\end{equation}

\subsection{Complex coordinates of the complexified tangent space}
On the complexified tangent space $\C T_{p}$ at a point $p\in M$,
we shall consider the complex coordinates
\begin{equation}\Label{}
dz=dx+idy\in\C^n_z, \, \bar dz=dx-idy\in\C^n_{\bar z}, \quad
dw\in \C^d_w \text{ with } du=\Re dw,
\end{equation}
where $\C^n_z$ and $\C^n_{\bar z}$
denote respectively the $(1,0)$ and $(0,1)$ spaces
with respect to the standard CR structure on $\C^n$.
Then our preliminary normalization \eqref{preli-norm}
is expressed by

\begin{equation}\Label{preli-norm1}
H^{1,0}_0 = \C^n_z\times\{0\}.
\end{equation}

Now consider a general point $p\in M$.
If $p$ is sufficiently close to $0$,
the subspace $H^{1,0}_p\subset \C T_p$ is the graph of a uniquely
determined complex-linear map
\begin{equation}\Label{mapL}
L(p)\colon \C^n_z\to \C^n_{\bar z}\times\C^d_{w},
\end{equation}
i.e.\  $H^{1,0}_p$ is given by $(d\bar z,dw)= L(p) dz$.

We shall also distinguish the components of $L$:
\begin{equation}\Label{}
L(p)=(L^{\bar z}(p),L^w(p)),
\quad L^{\bar z}(p)\colon \C^n_z \to \C^n_{\bar z},
\quad L^{w}(p)\colon \C^n_z \to \C^d_{w}.
\end{equation}
Then \eqref{preli-norm1} can be rewritten as
\begin{equation}\Label{preli-norm2}
L(0)=(L^{\bar z}(0),L^{w}(0))=0.
\end{equation}

\subsection{Evaluation along the Euler vector field}
In our normalization the following Euler (or radial) type vector field
will play an important role:
\begin{equation}\Label{euler}
e=e(z):=z\frac{\d}{\d z} = \sum_j z_j\frac{\d}{\d z_j} 
\in T^{1,0}_z\C^n.
\end{equation}
We consider $L(p)$ as a formal power series in $(z,\bar z,u)$ and evaluate it along $e$
to obtain a $\C^n_{\bar z}\times\C^d_w$-valued formal power series
\begin{equation}\Label{L-tilde}
\2L(z,\bar z,u):=L(z,\bar z,u) e(z) .
\end{equation}
We write $L_{z^{a}\bar z^b u^c}$ for the derivative at $0$
regarded as a multihomogeneous polynomial of degree $a$ in $z$,
$b$ in $\bar z$ and $c$ in $u$, which is given by
\begin{equation}\Label{multi-formula}
L_{z^{a}\bar z^b u^c}(z,\bar z,u):= \sum 
L_{z_{i_1}\ldots z_{i_a} \bar z_{j_1}\ldots \bar z_{j_b} u_{k_1} \ldots z_{k_c}}(0) 
z_{i_1}\cdots z_{i_a} \bar z_{j_1} \cdots \bar z_{j_b} u_{k_1} \cdots u_{k_c}
\end{equation}
in terms of the partial derivatives of $L$,
where the summation is taken over all collections of indices $i_1,\ldots,i_a$, $j_1,\ldots, j_b\in\{1,\ldots,n\}$, $k_1,\ldots, k_c\in\{1,\ldots,d\}$.
Then
\begin{equation}\Label{der-tilde}
\2L_{z^{a}\bar z^b u^c}(z,\bar z,u) = 
aL_{z^{a-1}\bar z^b u^c}(z,\bar z,u)e(z)
\end{equation}
for all integers $a,b,c\ge0$.
For the convenience of notation, 
we shall allow negative values of $a,b,c$
but always assume $h_{z^{a}\bar z^b u^c}=0$
for any function $h$ whenever any of $a,b,c$ is negative.

We shall say that two (formal) almost CR structures on $\C^n\times\R^d$
corresponding to formal maps $L,L'$,
{\em coincide along the Euler vector field},
if their evaluations along $e$ coincide, i.e.\ if $\2L\equiv\2L'$.

\subsection{Relation with the non-integrability tensor}
We next calculate the non-integrability tensor in terms of the above map $L$ defining
an almost CR-manifold
$(M,H,J)$, normalized as in \eqref{preli-norm2}.  At a point $p\in M$, it is given by the antisymmetric map
$\6N\colon H^{1,0}_p\times H^{1,0}_p\to \C T_p/ H^{1,0}_p$ in \eqref{N}
induced by the Lie brackets of $(1,0)$ vector fields.
Using the map \eqref{mapL}
we can choose $(1,0)$ vector fields of the form
$X=(\xi,L \xi)$ with $\xi$ being a constant vector field in $\C^n$.
Then for $X=(\xi,L\xi)$ and $Y=(\eta,L\eta)$, we have
\begin{equation}\Label{nij}
\begin{split}
\6N(X,Y) &= L_z(\xi;\eta) + L_{\bar z}(L^{\bar z}\xi;\eta)
+ L_u(L^{w}\xi;\eta)\\
&-L_z(\eta;\xi) - L_{\bar z}(L^{\bar z}\eta;\xi)
- L_u(L^{w}\eta;\xi),
\end{split}
\end{equation}
where we have adapted the notation
 \begin{equation}\Label{Lz}
L_z(\xi;\eta) = L_z(\xi)\eta
\end{equation}
for the derivative of $L$ in the direction of $\xi$ evaluated at $\eta$
and analogous notation for $L_{\bar z}$ and $L_u$,
and where $\C T_p/ H^{1,0}_p$ is  identified with $\C^n_{\bar z}\times\C_w^d$ in the obvious way.

The following lemma shows that when the given almost CR structures is integrable
(i.e.\ a CR structure), it is completely determined by the evaluation $\2L$ along the Euler vector field.
\bl\Label{integrable}
If two formal CR structures (i.e.\ integrable ones) of type $(n,d)$ on $M=\C^n\times\R^d$ are normalized as in \eqref{preli-norm1} and coincide along the Euler vector field \eqref{euler}, then they coincide as formal maps.
\el

\bpf
Denote by $L$ and $\3L$ the formal power series maps \eqref{mapL} corresponding to the given almost CR structures. 
We shall prove the coincidence of $L$ and $\3L$ at $0$ up to order $k$ by induction on $k$.
Suppose that all derivatives of $L$ and $\3L$ of order less than $k$ coincide at the origin.
Fix nonnegative integers $a,b,c$ with $a+b+c=k$ and consider the derivatives
$$D_{z^a\bar z^b u^c} L, D_{z^a\bar z^b u^c} \3L\colon (\C^n_z)^a\times(\C^n_{\bar z})^b\times(\R^d_u)^c\times\C^n_z\to \C^n_{\bar z}\times \C^d_w.$$
Then the coincidence of the given CR structures along the Euler vector field $e(z)$ implies\begin{equation}\Label{coinc}
D_{z^a\bar z^b u^c} L(z,\ldots,z,\bar z,\ldots,\bar z,u,\ldots,u;z)=
D_{z^a\bar z^b u^c} \3L(z,\ldots,z,\bar z,\ldots,\bar z,u,\ldots,u;z).
\end{equation}
In case $a=0$, this immediately implies $D_{\bar z^b u^c} L= D_{\bar z^b u^c} \3L$. Otherwise we apply the derivative $D_{z^{a-1}\bar z^b u^c}$ at the origin to \eqref{nij} written for $L$ and $\3L$. Since $L(0)=\3L(0)=0$ by the assumption, the corresponding derivative of the sum
$L_{\bar z}(L^{\bar z}\xi;\eta)
+ L_u(L^{w}\xi;\eta)-L_{\bar z}(L^{\bar z}\eta;\xi)
- L_u(L^{w}\eta;\xi)$
only involves nonzero contributions of derivatives of $L$ of order less than $k$
and similar property holds for $\3L$.
Thus, by the induction hypothesis, the derivatives of these sums are the same for $L$ and $\3L$. Now subtracting the derivatives of \eqref{nij} for $L$ and $\3L$ and using the assumption that  both CR structures are integrable, we conclude that the multilinear function
$D_{z^a\bar z^b u^c} L- D_{z^a\bar z^b u^c} \3L$ is invariant under exchanging
the first and the last arguments. Hence it is symmetric in all its $\C^n_z$-arguments
and therefore \eqref{coinc} implies $D_{z^a\bar z^b u^c} L= D_{z^a\bar z^b u^c} \3L$ as desired.
\epf

\subsection{Coordinate changes and basic identities}
We consider arbitrary formal power series coordinate changes
of the form
\begin{equation}\Label{fg}
z'=z+f(z,\bar z,u), \quad u'=u+g(z,\bar z,u),
\quad (f,g)=O(|(z,\bar z,u)|^2).
\end{equation}
We write $L'$ for the map \eqref{mapL}
corresponding to the new coordinates $(z',\1{z'},w')$.
Then $L$ and $L'$ are related by the following basic identities:

\begin{equation}\Label{basic}
\begin{split}
g_z + g_{\bar z} L^{\bar z} + (\id  + g_u) L^w &
= L'^{w'}(\id + f_z + f_{\bar z} L^{\bar z} + f_u L^w),\\
\bar f_z + (\id + \bar f_{\bar z}) L^{\bar z} + {\bar f}_u L^w &
= L'^{\bar z'}(\id + f_z + f_{\bar z} L^{\bar z} + f_u L^w).
\end{split}
\end{equation}
Here the derivatives of $f,\bar f,g$ and the maps $L^{\bar z}, L^w$
are taken at $(z,\bar z,u)$
and the maps $L'^{\bar z'},L'^{w'}$ at
$(z',\1{z'},u')=(z,\bar z,u)+(f,\bar f,g)(z,\bar z,u)$.

We shall use the evaluation of \eqref{basic} 
along the Euler vector field:
\begin{equation}\Label{basic-tilde}
\begin{split}
g_z e + g_{\bar z} \2L^{\bar z} + (\id  + g_u) \2L^w &
= L'^{w'}(e + f_z e + f_{\bar z} \2L^{\bar z} + f_u \2L^w),\\
\bar f_z e + (\id + \bar f_{\bar z}) \2L^{\bar z} + {\bar f}_u \2L^w &
= L'^{\bar z'}(e + f_z e + f_{\bar z} \2L^{\bar z} + f_u \2L^w).
\end{split}
\end{equation}

\section{Partial normal form for any CR-dimension and -codimension}
In the sequel all derivatives will be assumed evaluated at $0$ unless specified otherwise.
The derivatives of the maps $f$ and $g$ satisfy the following reality
conditions:
\begin{equation}\Label{reality}
\1{f_{z^a \bar z^b u^c}} = \bar f_{z^b \bar z^a u^c},
\quad \1{g_{z^a \bar z^b u^c}} = g_{z^b \bar z^a u^c}.
\end{equation}

\subsection{First order normalization}\Label{first}
We first take the derivatives of \eqref{basic} in $\bar z$,
evaluate at $0$ and use the vanishing
\eqref{preli-norm2} for both $L$ and
$L'$ as well as \eqref{fg} to obtain
\begin{equation}\Label{2-1}
g_{z\bar z} + L^w_{\bar z} = L'^{w'}_{\bar z'},
\quad \bar f_{z\bar z} + L^{\bar z}_{\bar z}
= L'^{\bar z'}_{\bar z'}.
\end{equation}
In the first equation $g_{z\bar z}$
is an arbitrary hermitian bilinear map in view of \eqref{reality}.
Hence we can use it to eliminate the hermitian part of
$L^w_{\bar z}(\bar\xi;\eta)$ (where we use the same notation as in \eqref{Lz}) and thus assume the latter to be antihermitian.
On the other hand, the term $\bar f_{z \bar z}$
in the second equation is completely arbitrary
and hence can be used to eliminate $L^{\bar z}_{\bar z}$
completely.

Similarly we take the derivatives of \eqref{basic} in $z$,
evaluate at $0$ and use the vanishing conditions
\eqref{preli-norm2} and \eqref{fg} to obtain
\begin{equation}\Label{2-2}
g_{z^2} + L^w_{z} = L'^{w'}_{z'},
\quad \bar f_{z^2} + L^{\bar z}_{z}
= L'^{\bar z'}_{z'}.
\end{equation}
This time both terms $g_{z^2}$ and $\bar f_{z^2}$
are arbitrary symmetric and can be used to eliminate
the symmetric parts of the bilinear forms $L^w_{z}(\xi;\eta)$ and 
$L^{\bar z}_{z}(\xi;\eta)$.
Hence we can normalize both forms to be antisymmetric.

Putting everything together, we obtain the normalization:
\begin{equation}\Label{norm1}
\begin{split}
L^{\bar z}_{\bar z}(\bar\xi;\eta) =0, \quad
&L^w_{\bar z}(\bar \xi;\eta) = - \1{L^w_{\bar z}(\bar \eta;\xi)},\\
L^{\bar z}_{z}(\xi;\eta) = -L^{\bar z}_{z}(\eta;\xi),\quad
&L^w_{z}(\xi;\eta) = -L^w_{z}(\eta;\xi).
\end{split}
\end{equation}
Note that in this normalization,
$\frac1{i} L^w_{\bar z}$ represents
the (hermitian) {\em Levi form} of the given
almost CR structure at $0$
and $2 L_{z}= (2L^{\bar z}_{z},2L^w_{z})$
the (antisymmetric) {\em non-integrability tensor} at $0$, i.e.\ 
\begin{equation}\Label{Levi-nonint}
\frac1{i} L^w_{\bar z} = \6L^{1,0}_0, \quad 2L_z = \6N_0.
\end{equation}
Furthermore, if both $L$ and $L'$ are normalized
as in \eqref{norm1},
it follows from \eqref{2-1} and \eqref{2-2} that
\begin{equation}\Label{vanish-low}
g_{z\bar z}=0, \quad g_{z^2}=0, \quad
\bar f_{z\bar z}=0, \quad \bar f_{z^2}=0.
\end{equation}

In view of \eqref{der-tilde},
the normalization \eqref{norm1} can also be rewritten 
in terms of $\2L$:
\begin{equation}\Label{}
\2L^{\bar z}_{z\bar z} =0, \quad \Re\2L^w_{z\bar z} =0,\quad
\2L_{z^2} =0.
\end{equation}

\subsection{Higher order expansion}\Label{high}
We now take arbitrary higher order derivatives of \eqref{basic-tilde}
that we regard as multi-homogeneous polynomials in $(z,\bar z,u)$ as in \eqref{multi-formula}.
As before, each derivative of $\2L$ is taken at $0$.

For every $a,b,c\ge 0$, we
differentiate \eqref{basic-tilde} $a$ times in $z$,
$b$ times in $\bar z$ and $c$ times in $u$
and evaluate at $0$.
In view of the vanishing in
 \eqref{preli-norm2} and \eqref{fg}, all terms in
\eqref{basic-tilde} (other than $\id$) vanish at $0$.
Hence we obtain:
\begin{equation}\Label{derived1}
\begin{split}
ag_{z^{a}\bar z^b u^c} + \2L^w_{z^{a}\bar z^b u^c} &=
\2L'^{w'}_{z'^{a}\bar z'^b u'^c} +
P_{a,b,c}(f_*,\bar f_*,g_*,\2L_*,L'^{w'}_{*}),\\
a\bar f_{z^{a}\bar z^b u^c} + \2L^{\bar z}_{z^{a}\bar z^b u^c} &=
\2L'^{\bar z'}_{z'^{a}\bar z'^b u'^c} + Q_{a,b,c}(f_*,\bar
f_*,g_*,\2L_{*},L'^{\bar z'}_{*}),
\end{split}
\end{equation}
where $P_{a,b,c}$ (resp.\ $Q_{a,b,c}$) is a polynomial in (the
components of) the derivatives (denoted by the subscript ``$*$") of $f$, $\bar f$, $g$ and $\2L$ 
of order less than $a+b+c$ 
and of $L'^{w'}$ 
(resp.\ $L'^{\bar z'}$) of order less than $a+b+c -1$.
(Note that the derivatives of $\2L$ of order $a+b+c$
are related to the derivatives of $L$ of order $a+b+c-1$
via \eqref{der-tilde}.)

Assuming all arguments of $Q_{a,b,c}$ being fixed, we can uniquely
choose $\bar f_{z^{a}\bar z^b u^c}$ in the second identity to make
$\2L'^{\bar z'}_{z'^{a}\bar z'^b u'^c}=0$ for $a\ge 1$,
whereas we already have $\2L'^{\bar z'}_{\bar z'^b u'^c}=0$
by \eqref{der-tilde} (recall that all derivatives are taken at $0$). 
This way we can uniquely
determine by induction all derivatives $\bar f_{z^{a}\bar z^b u^c}$
with $a\ge1$, i.e.\ all
derivatives of $\bar f$ except the pure ones 
$\bar f_{\bar z^b u^c}=\1{f_{z^bu^c}}$,
the latter being treated as free parameters.
Finally we use \eqref{basic} to determine all derivatives of $L'$
of order $a+b+c-1$ and complete the induction step.
Summarizing, we obtain:

\bl\Label{f-normal}
For every choice of the derivatives $g_{z^a\bar z^b u^c}$
and of pure derivatives $f_{z^a u^c}$, there
exists an unique choice of the remaining derivatives 
$f_{z^a \bar z^{b} u^c}$, $b\ge1$, such that
$\2L'^{\bar z'}_{z'^{a}\bar z'^b u'^c}=0$ for
all $a,b,c\ge 0$, i.e.\ $\2L'^{\bar z'}\equiv 0$.
\el


We shall next proceed similarly with the first identity in
\eqref{derived1}. In view of the reality conditions,
for each term $g_{z^{a}\bar z^b u^c}$ with $a\ne b$,
its conjugate appears in the other idenitity with $(a,b,c)$ replaced
by $(b,a,c)$. Thus we cannot eliminate both terms
$\2L'^{w'}_{z'^{a}\bar z'^b u'^c}$ and $\2L'^{w'}_{z'^{b}\bar
z'^{a} u'^c}$ simultaneously. However,
we can determine $g_{z^{a}\bar z^b u^c}$ uniquely
by eliminating the sum of the first one
and the conjugate of the second, i.e.\ by making
\begin{equation}\Label{conj1}
\2L'^{w'}_{z'^{a}\bar z'^b u'^c}+\1{\2L'^{w'}_{z'^{b}\bar z'^{a}
u'^c}}=0, \quad a\ne b.
\end{equation}
The sum in \eqref{conj1} is a certain derivative of
$\Re \2L'^{w'}(z',\bar z',u')$.

On the other hand, if $b=a>0$, both identities coincide but the
derivative $g_{z^a\bar z^a u^c}$ must be real. Hence it
is uniquely determined by eliminating the real part of
$\2L'^{w'}_{z'^{b}\bar z'^b u'^c}$:
\begin{equation}\Label{conj2}
\Re \2L'^{w'}_{z'^{a}\bar z'^a u'^c}=0, \quad a\ge1.
\end{equation}
Again, the real parts in \eqref{conj2} appear to be derivatives
of $\Re \2L'^{w'}(z',\bar z',u')$.


This way we determine all derivatives of $g$ except $g_{u^c}$.
The latters are to be treated as free parameters.
Finally we use \eqref{basic} to determine all derivatives of $L'$
of order $a+b+c-1$ and complete the induction step.
Summarizing, we obtain:
\bp\Label{partial}
For every formal power series $L(z,\bar z,u)\colon \C^n_z\to \C^n_{\bar z}\times\C^d_w$ 
without constant terms and every formal power series
$f_0(z,u)$ and $g_0(u)$ without constant and linear terms,
there exist unique formal power series $f(z,\bar z,u)$
and  $g(z,\bar z,u)$ without constant and linear terms
such that $f(z,0,u) = f_0(z,u)$ and $g(0,0,u)=g_0(u)$
and such that the map $\id+(f,g)$ transforms $L$ into $L'$
satisfying the normalization
\begin{equation}\Label{conj2'}
\2L'^{\bar z'}(z',\bar z',u')\equiv0, \quad
\Re \2L'^{w'}(z',\bar z',u')\equiv0.
\end{equation}
\ep

\section{Quasi CR embeddings}
\Label{quasi}
Recall that any real-analytic (integrable) CR structure
of CR dimension $n$ and CR codimension $d$ admits locally
a CR embedding into $\C^{n+d}$ inducing the given CR structure.
Vice versa, any almost CR structure induced by a CR embedding
into $\C^N$ (for arbitrary $N$) is automatically integrable.
Hence we clearly cannot have CR embeddings into any $\C^N$
for nonintegrable almost CR structures.
Here we propose a more general notion of {\em quasi CR embeddings}
that works for any almost CR structure and yields ``true" CR embeddings
whenever the almost CR structure is integrable.

\bd
A {\em quasi CR embedding} at a point $p$ 
of an almost CR structure of CR dimension $n$ and CR codimension $d$
 into $\C^n_z\times\C^d_w$ 
 is a (formal) embedding as a submanifold $M\subset \C^n_z\times\C^d_w$
with $p=0$, $T_0 M=\C^n\times\R^d$, 
such that the transformed almost CR structure on $M$
coincides along the Euler vector field $e$ (given by \eqref{euler})
with the one induced by the embedding,
where we regard $z$ and $u=\Re w$ 
as intrinsic coordinates on $M$.
(That is $\2L\equiv\2L'$ holds in the notation \eqref{L-tilde}
for the corresponding maps $L$ and $L'$ representing
the two almost CR structures on $M$.)
\ed

\subsection{Induced CR-structure for a graph}
Here we consider a formal (generic) submanifold $M$ in $\C^{n+d}$ given as
graph of a formal map $\phi\colon \C^n\times\R^d \to \R^d$, i.e.\
\begin{equation*}
M=\{(z,w)\in \C^n\times\C^d : \Im w = \phi(z,\bar z,\Re w)\}.
\end{equation*}
Then $(1,0)$ bundle $H^{1,0}M= 
T^{1,0}\C^{n+d} \cap \C TM$ is the annihilator of the (vector-valued) forms
\begin{equation}\Label{forms}
d\bar z, \quad \partial\Big(\frac{w-\bar w}{2i}
- \phi\big(z,\bar z, \frac{w+\bar w}{2}\big)\Big)=
\frac1{2i}((\id-i\phi_u)dw - 2i\phi_z dz).
\end{equation}
Hence, in terms of $(z,\bar z,u)$ with $u=\Re w$
regarded as coordinates on $M$,
the bundle $H^{1,0}$ is the annihilator of the pullbacks
of the forms \eqref{forms} under the map
$(z,u)\mapsto (z, u+i\phi(z,\bar z,u))$, i.e.\ it is the annihilator of
the forms
\begin{equation*}
d\bar z, \quad
(\id-i\phi_u)(\id+i\phi_u)du - (\id + i\phi_u) i\phi_z dz,
\end{equation*}
or, equivalently, since $\id-i\phi_u$ and $\id+i\phi_u$ commute, of the forms
\begin{equation}
d\bar z, \quad
(\id-i\phi_u)du - i\phi_z dz.
\end{equation}
Therefore the map $L$ is given by
\begin{equation}\Label{phi-L}
L = (L^{\bar z}, L^w)= \big(0,  i (\id -i\phi_u)^{-1}\phi_z\big).
\end{equation}

We now follow this construction backwards, i.e.\ begin with $L$
and reconstruct the function $\phi(z,\bar z,u)$.
In general, when the given almost CR structure is nonintegrable,
we cannot expect it to be induced by an embedding in a complex vector
space.
However we shall see that we can still find $\phi$ satisfying 
\eqref{phi-L} along the Euler vector field.   
In accordance with our normalization \eqref{preli-norm1},
we shall assume
\begin{equation}\Label{phi-vanish}
\phi(z,\bar z,u)=O(|(z,\bar z,u)|^2).
\end{equation}
Then evaluating \eqref{phi-L} along the Euler vector field 
\eqref{euler}
 and differentiating at $0$, we obtain
\begin{equation}\Label{der-lphi}
\2L^w_{z^{a}\bar z^b u^c} = i a \phi_{z^{a}\bar z^b u^c} + R_{a,b,c}(\phi_*)
\end{equation}
in the notation \eqref{L-tilde},
where $R_{a,b,c}$ is a polynomial in (the components of) 
the derivatives
$\phi_*$ of $\phi$ at $0$ of order less than $a+b+c$.
We shall consider functions $\phi$ satisfying 
\begin{equation}\Label{noharm}
\phi_{\bar z^b u^c}=0.
\end{equation}
This corresponds to the choice of normal coordinates
where the defining equation has no harmonic terms.
Given $\2L^w$, we use \eqref{der-lphi}
to obtain by induction on $a+b+c$ formulas for all
derivatives $\phi_{z^{a}\bar z^b u^c}$ other than those in \eqref{noharm}, i.e.\ with $a\ge1$:
\begin{equation}\Label{phi-recover}
\phi_{z^{a}\bar z^b u^c} = \frac1{ia}
\2L^w_{z^{a}\bar z^b u^c} +
S_{a,b,c}(\2L^w_{*}),
\end{equation}
where $S_{a,b,c}$ is a polynomial in the derivatives of $\2L^w$
of order less than $a+b+c$.

We now consider a transformation \eqref{fg} sending $L$ into $L'$
and look for a $\C^d$-valued function $\phi(z',\bar z',u')$
in the new coordinates providing a
local embedding of $\C^n\times\R^d$ into $\C^n\times\C^d$ given by
$(z',u')\mapsto (z',u'+i\phi(z',\bar z',u'))$, with desired properties.
Rewriting \eqref{phi-recover} for $L'$, we obtain, for $a\ge1$, 
\begin{equation}\Label{phi-recover'}
\phi_{z'^{a}\bar z'^b u'^c} = \frac1{ia}
\2L'^{w'}_{z'^{a}\bar z'^b u'^c} +
S_{a,b,c}(\2L'^{w'}_{*}).
\end{equation}
Solving the first equation in \eqref{derived1} for 
$\2L'^{w'}_{z'^{a}\bar z'^b u'^c}$ and substituting into 
\eqref{phi-recover'} we obtain by induction on $a+b+c$:
\begin{equation}\Label{phi-eq}
\phi_{z'^{a}\bar z'^b u'^c} = \frac1{i}g_{z^{a}\bar z^b u^c} 
+ T_{a,b,c}(f_*,\bar f_*,g_*,\2L'_{*}),
\end{equation}
where $T_{a,b,c}$ is a polynomial in the derivatives of $f$, $g$
of order less than $a+b+c$ and derivatives of $\2L'$
of order less than $a+b+c+1$.
Here we drop the dependence on $L$ which is assumed to be given and fixed.
We shall now use \eqref{phi-eq} along with the second
equation in \eqref{derived1} to determine uniquely
the derivatives $g_{z^{a}\bar z^b u^c}$ with $(a,b)\ne0$ and 
$\bar f_{z^{a}\bar z^b u^c}$ with $a\ne0$ as in \S\ref{high}
via the normalization conditions 
\begin{equation}\Label{normal-phi}
\2L'^{\bar z'}(z',\bar z',u')\equiv0, \quad \Im \phi(z',\bar z',u')\equiv0.
\end{equation}
As before, we complete the induction step 
by using \eqref{basic} to determine all derivatives of $L'$
of order $a+b+c-1$.

Summarizing and taking \eqref{noharm} into account, we obtain:
\bp\Label{quasi-prop}
For every formal power series 
$L(z,\bar z,u)\colon \C^n_z\to\C^n_{\bar z}\times\C^d_w$ without constant terms
and $f_0(z,u)\in\C^n$ and $g_0(u)\in\R^d$ without constant and linear terms,
there exist unique formal power series
$f(z,\bar z,u)\in\C^n$, $g(z,\bar z,u)\in\R^d$ and $\phi(z',\bar z',u')\in\R^d$  without constant and linear terms satisfying
$f(z,0,u) = f_0(z,u)$,  $g(0,0,u)=g_0(u)$ and $\phi(z',0,u')= 0$,
such that the almost CR structure given by $L$ admits 
the quasi CR embedding at $0$ as the submanifold $M'\subset\C^n_{z'}\times\C^d_{w'}$ given by
$$\Im w' = \phi(z',\bar z', \Re w')$$
via the map $(z,u)\mapsto (z+f(z,\bar z,u), u+g(z,\bar z,u))$, 
where $(z',u')$ are regarded as intrinsic coordinates on $M'$.
\ep

\br\Label{forms-different}
Note that even though the normalization conditions of Proposition~\ref{quasi-prop} look similar to \eqref{conj2'}, the two normalizations are different in general. For instance, consider the hypersurface given by $\Im w=\phi(z,\bar z,u)$ with $\phi(z,\bar z,u)=z\bar z+uz^4\bar z^4$. It is in the normal form of Proposition~\ref{quasi-prop} but 
$$\2L^w=\frac{i\phi_z e}{1-i\phi_u}=\frac{i(z\bar z+4uz^4\bar z^4)}{1-iz^4\bar z^4}$$
does not satisfy \eqref{conj2'}. 
\er

\section{Intrinsic normal form for hypersurface type almost CR structures}
\Label{hyper-form}
We now restrict our study the almost CR structures of hypersurface type,
i.e.\ those with CR codimension $1$.
\bd\Label{strongly}
We call an almost CR structure of hypersurface type
{\em strongly nondegenerate} if both its Levi form
and the linear combination $6i\6L+\6N^w$ are
nondegenerate, i.e.\ if, in addition to the nondegeneracy of the Levi
form  $\6L$, one has
\begin{equation}\Label{}
6i\6L(\bar\xi,\eta)+\6N^w(\xi,\eta) = 0\, \text{ for all } \eta
\quad\Rightarrow\quad \xi=0,
\end{equation}
where we assume the normalization \eqref{norm1}.
\ed
If the Levi form is positive definite, the almost CR structure
is always strongly nondegenerate.
Indeed, if $\xi$ is such that 
$6i\6L(\bar\xi,\eta) + \6N^w(\xi,\eta) = 0$
for all $\eta$, we have, in particular,
$6i\6L(\bar\xi,\xi) + \6N^w(\xi,\xi) = 0$.
Since $\6N$ is antisymmetric, $\6N(\xi,\xi)=0$ and therefore $\6L(\bar\xi,\xi)=0$.
In view of $\6L$ being positive definite, we obtain $\xi=0$ as desired.

\subsection{Trace decompositions}
We here recall trace decompositions
that play a fundamental role in \cite{CM}. 
Recall that the trace operator associated with the nondegenerate Levi
form $\frac1{i} \2L^w_{z\bar z}=\sum_{jk}c_{jk}d\bar z_j\otimes dz_k$
is given by $\tr:=\sum_{jk}c^{jk}\frac{\d^2}{\d \bar z_j\d z_k}$,
where $(c^{jk})$ is the inverse matrix of $(c_{jk})$.
(In particular, if $\2L^w_{z\bar z}=i\sum_j \eps_j d\bar z_j\otimes dz_j$ for $\eps=\pm1$, then $\tr=\sum_j \eps_j \frac{\d^2}{\d \bar z_j\d z_j}$.)
Then, for any integer $l\ge1$, any formal power series $p(z,\bar z,u)$ admits an unique decomposition
\begin{equation}\Label{decomp}
p(z,\bar z,u) = q(z,\bar z,u) (\2L^w_{z\bar z}(z,\bar z))^l + h(z,\bar z,u),
\end{equation}
where $q$ and $h$ are further power series and $h$ satisfies $\tr^l h\equiv0$.
Here $\tr^l$ stands for $\tr$ applied $l$ times.
Moreover, if $p$ is bihomogeneous of bidegree $(a,b)$ in $(z,\bar z)$,
then $h$ is also bihomogeneous of the same degree and
$q$ is bihomogeneous of degree $(a-l,b-l)$ (and is zero if $\min(a-l,b-l)<0$).
As in \cite{CM}, we shall use \eqref{decomp}
for $p$ of bidegrees $(2,2)$, $(3,2)$ and $(3,3)$ in $(z,\bar z)$
and $l$ equal $1$, $2$ and $3$ respectively.

\subsection{Weighted expansions}
As in \cite{CM} we assign weight $1$ to $z$, $\bar z$
and weight $2$ to $w$.
We shall assume the partial normalization as 
in Lemma~\ref{partial}.
For every $a,b,c\ge 0$, we
differentiate \eqref{basic-tilde} $a$ times in $z$,
$b$ times in $\bar z$ and $c$ times in $w$
and evaluate at $0$.
As before, we shall assume each derivative evaluated at $0$.
Then in view of the vanishing in
\eqref{preli-norm2} and \eqref{fg} we obtain
\begin{multline}\Label{basic1}
ag_{z^{a}\bar z^b u^c} + 
abg_{z^{a-1}\bar z^{b-1} u^{c+1}} \2L^w_{z\bar z}
+ \2L^w_{z^{a}\bar z^b u^c}
= \2L'^{w'}_{z'^{a}\bar z'^b u'^c} +
 P_{a,b,c}(f_*,\bar f_*,g_*,\2L_{*},L'^{w'}_{*})+\\
abL'^{w'}_{\bar z'}(\bar z;f_{z^{a}\bar z^{b-1} u^c}) +
aL'^{w'}_{\bar z'}(\bar f_{z^{a-1}\bar z^{b} u^c};z) +
a(a-2)L'^{w'}_{z'}(z;f_{z^{a-1}\bar z^{b} u^c})+\\
ab(b-1)L'^{w'}_{\bar z'}(\bar z;f_{z^{a-1}\bar z^{b-2} u^{c+1}})
\2L^w_{z\bar z} +
ab(a-1)L'^{w'}_{z'}(z;f_{z^{a-2}\bar z^{b-1} u^{c+1}})
\2L^w_{z\bar z},
\end{multline}
\begin{equation}\Label{basic2}
a\bar f_{z^{a}\bar z^b u^c} +
ab\bar f_{z^{a-1}\bar z^{b-1}u^{c+1}}\2L^w_{z\bar z}+
\2L^{\bar z}_{z^a\bar z^b u^c}
= \2L'^{\bar z'}_{z'^a\bar z'^b u'^c} +
Q_{a,b,c}(f_*,\bar f_*,g_*,\2L_{*},L'^{\bar z'}_{*}),
\end{equation}
where $P_{a,b,c}$ (resp.\ $Q_{a,b,c}$) is a polynomial in (the
components of) the derivatives of $f$, $\bar f$ and $L'^{w'}$ 
of weight less than $a+b+2c-1$ 
and of $g$ and $\2L$
of weight less than $a+b+2c$  
(resp. of\ $L'^{\bar z'}$ of weight less than $a+b+2c-1$ and $f$, $\bar f$, $g$ and $\2L$
of weight less than $a+b+2c$).
Here we have used the antisymmetry of $L_{z}(\cdot;\cdot)$ and 
$L'_{z'}(\cdot;\cdot)$  implying the cancellations
$$\2L_{z^2}=2L_{z}(z;z)=0, \quad L'_{z'}(f_{z^{a-1}\bar z^{b} u^c};z) +
L'_{z'}(z;f_{z^{a-1}\bar z^{b} u^c})=0.$$
We have also used \eqref{vanish-low} (to avoid terms like
$L'^{w'}_{z'^{a-1}\bar z'^{b-1} u'^{c+1}}g_{z\bar z}$ etc.).
Also recall that any derivative $h_{z^{a}\bar z^b u^c}$
with either of $a,b,c$ negative is assumed to be zero.

As before, assuming the arguments of $Q_{a,b,c}$ in
\eqref{basic2} being given,
the non-pure derivatives $\bar f_{z^{a} \bar z^b u^c}$
(i.e.\ those with $a\ge1$)
are uniquely determined by the first identity 
of our partial normalization
\eqref{conj2'},
i.e.\ from
\begin{equation}\Label{basic21}
\bar f_{z^{a}\bar z^b u^c} =
-b\bar f_{z^{a-1}\bar z^{b-1}u^{c+1}}\2L^w_{z\bar z}
-\frac{1}{a}\2L^{\bar z}_{z^a\bar z^b u^c}
+\frac1{a}Q_{a,b,c}(f_*,\bar f_*,g_*,\2L_{*},L'^{\bar z'}_{*}).
\end{equation}
However, this time we no more regard the pure derivatives $f_{z^a u^c}$
as free parameters but rather want to determine
them from the identity \eqref{basic1}
by adding further normalization conditions.
We shall assume $L$ being given and determine the derivatives of $g$ of a fixed weight $k$ and of $f$ of weight $k-1$. Then we use \eqref{basic} to determine all derivatives of $L'^{w'}$ of weight $k-1$ and of $L'^{\bar z'}$ of weight $k-2$. Thus our main inductive hypothesis for an integer $k\ge3$,
will be that we have already determined 
all derivatives of $g$ of weight
less than $k$ and all derivatives of $f$ and $L'^{w'}$
of weight less than $k-1$ and of $L'^{\bar z'}$ of weight less than $k-2$.
Our goal is then to determine the corresponding derivatives
of the next following weights.

\subsection{The diagonal terms normalization}\Label{diagonal}
We begin by considering the identities \eqref{basic1}
corresponding to $(a,b,c)$ equal $(1,1,s+1)$, $(2,2,s)$ 
and $(3,3,s-1)$.
Here $s\ge 0$ is such that $4+2s=k$ (i.e.\ the weight $a+b+2c$ is $k$).
(In case $s=0$ we only have the first two identities
and regard the third idenity as void.) We obtain
\begin{multline}\Label{basic1'}
g_{z\bar z u^{s+1}} + g_{u^{s+2}} \2L^w_{z\bar z}
+ \2L^w_{z\bar z u^{s+1}}
= \2L'^{w'}_{z'\bar z' u'^{s+1}} +
 P_{1,1,s+1}(f_*,\bar f_*,g_*,\2L_{*},L'^{w'}_{*})+\\
L'^{w'}_{\bar z'}(\bar z;f_{z u^{s+1}}) +
L'^{w'}_{\bar z'}(\bar f_{\bar z u^{s+1}};z)-
L'^{w'}_{z'}(z;f_{\bar z u^{s+1}}),
\end{multline}
\begin{multline}\Label{basic1''}
2g_{z^2\bar z^2 u^s} + 4g_{z\bar z u^{s+1}} L^w_{z\bar z}
+ \2L^w_{z^2\bar z^2 u^s}
= \2L'^{w'}_{z'^2\bar z'^2 u'^s} +
 P_{2,2,s}(f_*,\bar f_*,g_*,\2L_{*},L'^{w'}_{*})+\\
4L'^{w'}_{\bar z'}(\bar z;f_{z^2\bar z u^s}) +
2L'^{w'}_{\bar z'}(\bar f_{z\bar z^2 u^s};z) +
4L'^{w'}_{\bar z'}(\bar z;f_{z u^{s+1}})\2L^w_{z\bar z} +
4L'^{w'}_{z'}(z;f_{\bar z u^{s+1}})\2L^w_{z\bar z},
\end{multline}
\begin{multline}\Label{basic1'''}
3g_{z^3\bar z^3 u^{s-1}} + 9g_{z^2\bar z^2 u^s} \2L^w_{z\bar z}
+ \2L^w_{z^3\bar z^3 u^{s-1}}
= \2L'^{w'}_{z'^3\bar z'^3 u'^{s-1}} +
 P_{3,3,s-1}(f_*,\bar f_*,g_*,\2L_*,L'^{w'}_{*})+\\
9L'^{w'}_{\bar z'}(\bar z;f_{z^{3}\bar z^{2} u^{s-1}}) +
3L'^{w'}_{\bar z'}(\bar f_{z^2\bar z^{3} u^{s-1}};z) +
3L'^{w'}_{z'}(z;f_{z^2\bar z^{3} u^{s-1}}) +
\\
18L'^{w'}_{\bar z'}(\bar z;f_{z^2\bar z u^{s}})\2L^w_{z\bar z} +
18L'^{w'}_{z'}(z;f_{z\bar z^2 u^{s}})\2L^w_{z\bar z}.
\end{multline}

The first identity has the non-pure derivative $f_{\bar z u^{s+1}}$
that we can substitute from \eqref{basic21},
where all unknown terms on the right-hand side are of lower weights and
hence are already determined.
We now impose the condition
\begin{equation}\Label{n1}
\2L'^{w'}_{z'\bar z' u'^{s+1}}=0
\end{equation}
(rather than requiring to vanish only the real part as in 
\eqref{conj2'}).
Taking real and imaginary parts of \eqref{basic1'}
and using reality \eqref{reality} as well as the facts
that $\Re\2L^w_{z\bar z} =0$ and 
that $L^{w'}_{\bar z'}(\cdot;\cdot)$
is nondegenerate and antihermitian (see \eqref{norm1}), we see that 
by choosing $g_{z\bar zu^{s+1}}$ and $f_{z u^{s+1}}$ suitably, \eqref{n1}
can always be achieved and uniquely determines
$g_{z\bar zu^{s+1}}$ as well as the expression
\begin{equation}\Label{e1}
E_1:=g_{u^{s+2}} \2L^w_{z\bar z} -
2i\Im L'^{w'}_{\bar z'}(\bar z;f_{z u^{s+1}}).
\end{equation}

The next identity \eqref{basic1''} has the non-pure derivatives 
$f_{z^2\bar z u^s}$, $\bar f_{z\bar z^2 u^s}$
and $f_{\bar z u^{s+1}}$ that we can substitute from \eqref{basic21}.
After the substitution, the identity becomes 
\begin{equation}\Label{basic11''}
\begin{split}
2g_{z^2\bar z^2 u^s} & + 4g_{z\bar z u^{s+1}} L^w_{z\bar z}
	+ \2L^w_{z^2\bar z^2 u^s} 
	\\
= & \2L'^{w'}_{z'^2\bar z'^2 u'^s} +
	(12L'^{w'}_{\bar z'}(\bar z;f_{z u^{s+1}}) -4L'^{w'}_{\bar z'}(\bar f_{\bar z u^{s+1}};z) )\2L^w_{z\bar z} +
\ldots,
\end{split}
\end{equation}
where the dots stand for the terms already determined.
The derivative $g_{z\bar zu^{s+1}}$ has also been determined
from the previous identity.
Taking the imaginary parts of both sides in \eqref{basic11''}
we see that $\Im \2L'^{w'}_{z'^2\bar z'^2 u'^s}$ is determined
up to addition of a multiple of
$$\Im \big(L'^{w'}_{\bar z'}(\bar z;f_{z u^{s+1}})\2L^w_{z\bar z}\big) = 
-i \big(\Re L'^{w'}_{\bar z'}(\bar z;f_{z u^{s+1}})\big)\2L^w_{z\bar z},$$
as $\2L^w_{z\bar z}$ is purely imaginary.
Since the Levi form  
$\frac1{i}L'^{w'}_{\bar z'}(\cdot;\cdot)$ is nondegenerate,
we can always achieve the normalization
\begin{equation*}
\tr (\Im \2L'^{w'}_{z'^2\bar z'^2 u'^s}) = 0,
\end{equation*}
which determines uniquely the expression
\begin{equation}\Label{e2}
E_2:=\Re L'^{w'}_{\bar z'}(\bar z;f_{z u^{s+1}}).
\end{equation}
Here we have used the fact that both real and imaginary parts
of $L'^{w'}_{\bar z'}(\bar z;f_{z u^{s+1}})$ can be arbitrarily prescribed independently of each other by suitably choosing $f_{z u^{s+1}}$.

Next, suitably choosing $g_{z^2\bar z^2 u^s}$ in the real part of \eqref{basic11''}, 
we can also achieve
\begin{equation*}
\Re  \2L'^{w'}_{z'^2\bar z'^2 u'^s} =0.
\end{equation*}
The latter condition determines uniquely the expression
\begin{equation}\Label{e3'}
\begin{split}
2E_3: & = 2g_{z^2\bar z^2 u^s} -
	8\Re (L'^{w'}_{\bar z'}(\bar z;f_{z u^{s+1}})\2L^w_{z\bar z}) 
\\
& = 2g_{z^2\bar z^2 u^s} -
	(8i\Im L'^{w'}_{\bar z'}(\bar z;f_{z u^{s+1}}))\2L^w_{z\bar z}.
\end{split}
\end{equation}
Solving for $g_{z^2\bar z^2 u^s}$ and using \eqref{e1} we obtain
\begin{equation}\Label{g22}
g_{z^2\bar z^2 u^s} = g_{u^{s+2}} (\2L^w_{z\bar z})^2
-E_1 \2L^w_{z\bar z} + E_3.
\end{equation}
Then, in case $s=0$, the derivative $g_{u^2}$
determines $g_{z^2\bar z^2 u^s}$
and via \eqref{e1} and \eqref{e2}
both real and imaginary parts of 
$L'^{w'}_{\bar z'}(\bar z;f_{z u^{s+1}})$.
Since $L'^{w'}_{\bar z'}(\cdot;\cdot)$ is assumed nondegenerate,
the latters uniquely determine $f_{z u^{s+1}}$.
Thus we have determined all derivatives involved
except $g_{u^2}$, which is treated as free parameter.

Finally, in case $s\ge1$, we analyse the last identity \eqref{basic1'''}. 
Using \eqref{basic21} as before to substitute for all non-pure derivatives of $f$,
we obtain
\begin{multline}\Label{basic21'''}
3g_{z^3\bar z^3 u^{s-1}} + 9g_{z^2\bar z^2 u^s} \2L^w_{z\bar z}
+ \2L^w_{z^3\bar z^3 u^{s-1}}
= \2L'^{w'}_{z'^3\bar z'^3 u'^{s-1}} +\\
(90L'^{w'}_{\bar z'}(\bar z;f_{zu^{s+1}}) +
18L'^{w'}_{\bar z'}(\bar f_{\bar z u^{s+1}};z))(\2L^w_{z\bar z})^2 +\ldots.
\end{multline}
Taking real parts of both sides and choosing $g_{z^3\bar z^3 u^{s-1}}$ suitably, we see that
the condition
\begin{equation}\Label{}
\Re \2L'^{w'}_{z'^3\bar z'^3 u'^{s-1}} =0
\end{equation}
can be fulfilled and determines uniquely the expression
\begin{equation}
E_4:=3g_{z^3\bar z^3 u^{s-1}}-72\Re L'^{w'}_{\bar z'}(\bar z;f_{zu^{s+1}}) (\2L^w_{z\bar z})^2.
\end{equation}
Furthermore, since \eqref{e2} is determined, we conclude that 
we have determined $g_{z^3\bar z^3 u^{s-1}}$.

We now take the imaginary parts of both sides of \eqref{basic21'''},
where we substitute $g_{z^2\bar z^2 u^s}$ from \eqref{e3'}
and $\Im L'^{w'}_{\bar z'}(\bar z;f_{zu^{s+1}})$ from \eqref{e1}:
\begin{equation}\Label{basic22'''}
\Im\2L^w_{z^3\bar z^3 u^{s-1}}
= \Im\2L'^{w'}_{z'^3\bar z'^3 u'^{s-1}} -36i g_{u^{s+2}} (\2L^w_{z\bar z})^3
 +\ldots.
\end{equation}
We see that $\Im \2L'^{w'}_{z'^3\bar z'^3 u'^{s-1}}$
is determined up to a multiple  of $ig_{u^{s+2}} (\2L^w_{z\bar z})^3$.
Hence, choosing suitably $g_{u^{s+2}}$, we can always achieve
\begin{equation}\Label{}
\tr^3 (\Im \2L'^{w'}_{z'^3\bar z'^3 u'^{s-1}}) =0
\end{equation}
that determines $g_{u^{s+2}}$ uniquely.
As in case $s=0$ above, we see
that $g_{u^{s+2}}$ in turn determines $g_{z^2\bar z^2 u^s}$
via \eqref{g22} and $f_{z u^{s+1}}$ via \eqref{e1} and \eqref{e2},
where the nondegeneracy of $L'^{w'}_{\bar z'}(\cdot;\cdot)$ has been used.

We finally consider the remaining diagonal term expansions,
i.e.\ the identities \eqref{basic1} for $(a,b,c)$ with $a=b\ge4$. 
All derivatives of $f$  involved there have either been already determined or 
can be determined using \eqref{basic21}.
We then obtain the identities of the form
\begin{equation}\Label{age4}
ag_{z^{a}\bar z^{a} u^c} + a^2g_{z^{a-1}\bar z^{a-1} u^{c+1}} 
\2L^w_{z\bar z}
+ \2L^w_{z^a\bar z^{a} u^c}
= \2L'^{w'}_{z'^a\bar z'^{a} u'^c} + \ldots,
\end{equation}
where as before the dots stand for the terms that have 
been previously determined.
We can argue by induction on $a$ to see that
there is an unique choice of
the derivatives $g_{z^{a}\bar z^{a} u^c}$ such that
\begin{equation}\Label{}
\Re \2L'^{w'}_{z'^a\bar z'^{a} u'^c} =0.
\end{equation}
 
Summarizing, we have uniquely determined the derivatives
 $f_{z^a\bar z^{b}u^c}$ with $|a-b|=1$ of weight $k-1$ as well as
$g_{z^a\bar z^a u^c}$ of weight $k$
except the free parameter $g_{u^2}$, by the conditions
$\2L'^{\bar z'}\equiv0$ and
\begin{equation}\Label{n100}
\Re \2L'^{w'}_{z'^a\bar z'^{a} u'^c} =0,\quad
\Im \2L'^{w'}_{z'\bar z' u'^{c+1}} =0, \quad
\tr(\Im \2L'^{w'}_{z'^2\bar z'^2 u'^c}) =0, \quad
\tr^3(\Im \2L'^{w'}_{z'^3\bar z'^3 u'^{c}}) =0.
\end{equation}

\subsection{Normalization of the terms next to the diagonal}\Label{off}
Here we are going to normalize the terms of \eqref{basic1}
that correspond to $|a-b|=1$.
We begin by considering the identities for $(a,b,c)$
equal to $(1,0,s+1)$, $(2,1,s)$, $(1,2,s)$, $(3,2,s-1)$ and $(2,3,s-1)$
with $s\ge 0$ and $2s+3=k$ (i.e.\ the weight $a+b+2c=k$ as in
the previous section). In case $s=0$ the last two idenitities are
not present. Using \eqref{basic21} to eliminate all non-pure derivatives of $f$
as before in \S\ref{diagonal}, 
we obtain the identities
\begin{equation}\Label{b}
\begin{split}
g_{zu^{s+1}} + \2L^w_{z u^{s+1}}
=& \2L'^{w'}_{z' u'^{s+1}} +
L'^{w'}_{\bar z'}(\bar f_{u^{s+1}};z) - L'^{w'}_{z'}(z; f_{u^{s+1}})+
\ldots,\\
2g_{z^{2}\bar z u^s} + 2g_{z u^{s+1}} \2L^w_{z\bar z}
+ \2L^w_{z^2\bar z u^s}
=& \2L'^{w'}_{z'^2\bar z' u'^s} +
2L'^{w'}_{\bar z'}(\bar z;f_{z^{2} u^s}) \\&+
(- 2L'^{w'}_{\bar z'}(\bar f_{ u^{s+1}};z)  +
2L'^{w'}_{z'}(z;f_{u^{s+1}}))\2L^w_{z\bar z} + \ldots,\\
g_{z\bar z^2 u^s} + 2g_{\bar z u^{s+1}} \2L^w_{z\bar z}
+ \2L^w_{z\bar z^2 u^s}
=& \2L'^{w'}_{z'\bar z'^2 u'^s} +
L'^{w'}_{\bar z'}(\bar f_{\bar z^{2} u^s};z) +
4L'^{w'}_{\bar z'}(\bar z;f_{ u^{s+1}})\2L^w_{z\bar z}+\ldots,\\
3g_{z^{3}\bar z^2 u^{s-1}} + 6g_{z^2\bar z u^{s}} \2L^w_{z\bar z}
+ \2L^w_{z^3\bar z^2 u^{s-1}}
=& \2L'^{w'}_{z'^3\bar z'^2 u'^{s-1}} +
24L'^{w'}_{\bar z'}(\bar z;f_{z^2 u^{s}})\2L^w_{z\bar z} \\&+ 
(6L'^{w'}_{\bar z'}(\bar f_{ u^{s+1}};z)+
18L'^{w'}_{z'}(z;f_{u^{s+1}}))(\2L^w_{z\bar z})^2+\ldots,\\
2g_{z^{2}\bar z^3 u^{s-1}} + 6g_{z\bar z^2 u^s} \2L^w_{z\bar z}
+ \2L^w_{z^2\bar z^3 u^{s-1}}
=& \2L'^{w'}_{z'^2\bar z'^3 u'^{s-1}} -6L'^{w'}_{\bar z'}(\bar f_{\bar z^{2} u^s};z)
\2L^w_{z\bar z} \\&+ 
24L'^{w'}_{\bar z'}(\bar z;f_{ u^{s+1}}) (\2L^w_{z\bar z})^2 +\ldots,
\end{split}
\end{equation}
where as before the dots stand for the terms that have been already determined.

Our next normalization condition is
\begin{equation}\Label{n11}
\2L'^{w'}_{z' u'^{s+1}} = 0,
\end{equation}
which can always be achieved in view of the first identity in \eqref{b}
by suitably choosing $g_{zu^{s+1}}$,
and determines uniquely the quantity
\begin{equation}\Label{e11}
E_1:=  g_{zu^{s+1}} - L'^{w'}_{\bar z'}(\bar f_{u^{s+1}};z) +
 L'^{w'}_{z'}(z; f_{u^{s+1}}).
\end{equation}
Since $\2L'^{w'}(0,\bar z,u)\equiv0$ (following from the definition), 
\eqref{n11} is equivalent to
\begin{equation}\Label{n11'}
\2L'^{w'}_{z' u'^{s+1}}+\1{\2L'^{w'}_{\bar z' u'^{s+1}}} = 0,
\end{equation}
where the second term vanishes.

We next consider the second identity in \eqref{b} 
and twice the conjugate of the third.
Using the reality properties \eqref{reality} of $g$
and the antihermitian properties of $L'^{w'}_{\bar z'}(\cdot;\cdot)$ and $\2L^w_{z\bar z}$,
we obtain
\begin{multline}\Label{i1}
4g_{z^2\bar z u^s} -2 g_{z u^{s+1}} \2L^w_{z\bar z} + 
(\2L^w_{z^2\bar z u^s} + 2\1{\2L^w_{z\bar z^2 u^s}})=\\
 (\2L'^{w'}_{z'^2\bar z' u'^s} +  2\1{\2L'^{w'}_{z'\bar z'^2 u'^s}}) +
(6L'^{w'}_{\bar z'}(\bar f_{u^{s+1}};z)
+2L'^{w'}_{z'}(z;f_{u^{s+1}}))\2L^w_{z\bar z}+ \ldots.
\end{multline}
Similarly, substracting twice the conjugate of the third identity 
in \eqref{b} from the second, we obtain
\begin{multline}\Label{i2}
6g_{zu^{s+1}}\2L^w_{z\bar z} + (\2L^w_{z^2\bar z u^s}
- 2\1{\2L^w_{z\bar z^2 u^s}})=\\
( \2L'^{w'}_{z'^2\bar z' u'^s} -  2\1{\2L'^{w'}_{z'\bar z'^2 u'^s}}) +
4L'^{w'}_{\bar z'}(\bar z;f_{z^{2} u^s}) +
( -10L'^{w'}_{\bar z'}(\bar f_{u^{s+1}};z)+2L'^{w'}_{z'}(z;f_{u^{s+1}}))
\2L^w_{z\bar z}+ \ldots.
\end{multline}
Using the freedom in the choice of $g_{z^2\bar zu^s}$
in \eqref{i1} and of $f_{z^2u^s}$ in \eqref{i2}
together with nondegeneracy of $L'^{w'}_{\bar z'}(\cdot;\cdot)$
we can achieve the normalization
\begin{equation}\Label{n12}
\2L'^{w'}_{z'^2\bar z' u'^s} 
\pm  2\1{\2L'^{w'}_{z'\bar z'^2 u'^s}} =0,
\end{equation}
which determines uniquely the expressions
\begin{equation}\Label{e12}
\begin{split}
E_2 &:= 4g_{z^2\bar z u^s} - 2 g_{z u^{s+1}} \2L^w_{z\bar z}-
 (6L'^{w'}_{\bar z'}(\bar f_{u^{s+1}};z) +2L'^{w'}_{ z'}(z;f_{u^{s+1}}))\2L^w_{z\bar z},\\
E_3 &:= 6g_{z u^{s+1}}\2L^w_{z\bar z}
-4L'^{w'}_{\bar z'}(\bar z;f_{z^{2} u^s})
+(10L'^{w'}_{\bar z'}(\bar f_{u^{s+1}};z)- 2L'^{w'}_{ z'}(z;f_{u^{s+1}}))\2L^w_{z\bar z},
\end{split}
\end{equation}
Notice that in case $s=0$, we have $f_{u^{s+1}}=f_u=0$
in view of \eqref{fg}.
Then $g_{zu^{s+1}}$ is uniquely determined from \eqref{e11}
and both derivatives $g_{z^2\bar z u^s}$ and $f_{z^{2} u^s}$
are uniquely determined from respectively the first and the second identity \eqref{e12}.

In case $s\ge 1$ we have the last two identities in \eqref{b}.
We consider the sum of the first of them and the conjugate 
of the second as well as the difference between $2$ times the first and $3$ times the conjugate of the second:
\begin{equation}\Label{red}
\begin{split}
	5g_{z^3\bar z^2 u^{s-1}}& +  (\2L^w_{z^3\bar z^2 u^{s-1}} + \1{\2L^w_{z^2\bar z^3 u^{s-1}}}) \\
		=(&\2L'^{w'}_{z'^3\bar z'^2 u'^{s-1}} +  \1{\2L'^{w'}_{z'^2\bar z'^3 u'^{s-1}}})
		+ 18L'^{w'}_{\bar z'}(\bar z;f_{ z^2 u^s})\2L^w_{z\bar z} \\
	&+ (-18L'^{w'}_{\bar z'}(\bar f_{ u^{s+1}};z)+
		18L'^{w'}_{z'}(z;f_{u^{s+1}}))(\2L^w_{z\bar z})^2 +\ldots,\\
	30g_{z^2\bar z u^s}&\2L^w_{z\bar z}  +  (2\2L^w_{z^3\bar z^2 u^{s-1}} -
			3\1{\2L^w_{z^2\bar z^3 u^{s-1}}}) \\
		=(&2\2L'^{w'}_{z'^3\bar z'^2 u'^{s-1}} -
3\1{\2L'^{w'}_{z'^2\bar z'^3 u'^{s-1}}})
+ 66L'^{w'}_{\bar z'}(\bar z;f_{z^2 u^s})\2L^w_{z\bar z}
\\&+(12\cdot 7 L'^{w'}_{\bar z'}(\bar f_{ u^{s+1}};z)
+36L'^{w'}_{z'}(z;f_{u^{s+1}}))(\2L^w_{z\bar z})^2
+\ldots,
\end{split}
\end{equation}
Using the freedom in the choice of $g_{z^3\bar z^2 u^{s-1}}$
we can achieve
\begin{equation}\Label{n13}
\2L^w_{z^3\bar z^2 u^{s-1}} + \1{\2L^w_{z^2\bar z^3 u^{s-1}}} =0,
\end{equation}
which uniquely determines the expression
\begin{equation}\Label{e13}
E_4:= 5g_{z^3\bar z^2 u^{s-1}}
- 18L'^{w'}_{\bar z'}(\bar z;f_{ z^2 u^s})\2L^w_{z\bar z} +
(18L'^{w'}_{\bar z'}(\bar f_{ u^{s+1}};z)-
18L'^{w'}_{z'}(z;f_{u^{s+1}}))(\2L^w_{z\bar z})^2.
\end{equation}

The normalization of the difference
$(\2L'^{w'}_{z'^3\bar z'^2 u'^{s-1}} -
\1{\2L'^{w'}_{z'^2\bar z'^3 u'^{s-1}}})$ is trickier.
Solving \eqref{e11} for $g_{zu^{s+1}}$ and
the identities in \eqref{e12} for $g_{z^2\bar z u^s}$
and $L'^{w'}_{\bar z'}(\bar z;f_{z^{2} u^s})$ respectively
and substituting into the second identity in \eqref{red}, we obtain
\begin{equation}\Label{2-3}
\begin{split}
	-96&\big(3L'^{w'}_{\bar z'}(\bar f_{u^{s+1}};z)
		+L'^{w'}_{ z'}(f_{u^{s+1}};z)
		\big)(\2L^w_{z\bar z})^2
		+
		(2\2L^{w}_{z^3\bar z^2 u^{s-1}} -
		3\1{\2L^{w}_{z^2\bar z^3 u^{s-1}}}) \\
	&= (2\2L'^{w'}_{z'^3\bar z'^2 u'^{s-1}} -
3\1{\2L'^{w'}_{z'^2\bar z'^3 u'^{s-1}}})
+\ldots,
\end{split}
\end{equation}
where we have used the antisymmetry of $L'^{w'}_{ z'}(\cdot;\cdot)$.
Now, if the given almost CR structure is {\em strongly nondegenerate} (see Definition~\ref{strongly} and relation \eqref{Levi-nonint}),
we can choose uniquely $f_{u^{s+1}}$ to satisfy the normalization
condition
\begin{equation*}
\tr^2 (2\2L'^{w'}_{z'^3\bar z'^2 u'^{s-1}} -
3\1{\2L'^{w'}_{z'^2\bar z'^3 u'^{s-1}}})=0,
\end{equation*}
which in view of \eqref{n13} is equivalent to
\begin{equation}\Label{tr2'}
\tr^2 (\2L'^{w'}_{z'^3\bar z'^2 u'^{s-1}} -
\1{\2L'^{w'}_{z'^2\bar z'^3 u'^{s-1}}})=0.
\end{equation}
Once $f_{u^{s+1}}$ is determined, $g_{zu^{s+1}}$ is determined
via \eqref{e11}, $g_{z^2\bar z u^s}$ and $f_{z^2u^s}$ via \eqref{e12}
and subsequently $g_{z^3\bar z^2 u^{s-1}}$ via \eqref{e13}.

The remaining terms $g_{z^{a+1}\bar z^a u^c}$ 
are normalized inductively
in the same way as in the previous subsection by the conditions
\begin{equation}\Label{re}
\2L'^{w'}_{z'^{a+1}\bar z'^a u'^c} +
\1{\2L'^{w'}_{z'^{a}\bar z'^{a+1} u'^c}}=0.
\end{equation}
Summarizing we have obtained the normalization conditions
\begin{equation}\Label{n101}
\begin{split}
	\2L'^{w'}_{z'^{a+1}\bar z'^{a} u'^c} +&
		\1{\2L'^{w'}_{z'^{a}\bar z'^{a+1} u'^c}}=0,\quad
		\2L'^{w'}_{z'^2\bar z' u'^{c}} -
		\1{\2L'^{w'}_{z'\bar z'^2 u'^{c}}} =0, \\
	&\tr^2 (\2L'^{w'}_{z'^3\bar z'^2 u'^{c}} -
		\1{\2L'^{w'}_{z'^2\bar z'^3 u'^{c}}})=0, 
\end{split}
\end{equation}
that uniquely determine the derivatives $f_{u^c}$ and $f_{z^2u^c}$ of weight $k-1$
as well as $g_{z^{a+1} \bar z^a u^c}$ of weight $k$.

\subsection{Normalization of the remaining terms}
We now consider the identities \eqref{basic1}
that haven't been previously used, i.e.\ those corresponding
to $(a,b,c)$ with $|a-b|\ge2$ and $a+b+2c=k$.
We begin with the identitites for $(a,b,c)$ equal to
$(a,0,c)$, $(a+1,1,c-1)$, $(1,a+1,c-1)$ (with $a\ge2$),
where we use \eqref{basic21} as before to eliminate the non-pure
derivatives of $f$:
\begin{equation}\Label{bb0}
\begin{split}
ag_{z^{a}u^c}  +  \2L^w_{z^a u^c} &
	= \2L'^{w'}_{z'^a u'^c} + a(a-2)L'^{w'}_{z'}(z;f_{z^{a-1} u^{c}})+ \ldots,
	\\
(a+ & 1)g_{z^{a+1}\bar z u^{c-1}} + (a+1)g_{z^{a}u^{c}} \2L^w_{z\bar z}
	+ \2L^w_{z^{a+1}\bar z u^{c-1}} 
	\\
= & \2L'^{w'}_{z'^{a+1}\bar z' u'^{c-1}} +
	(a+1)L'^{w'}_{\bar z'}(\bar z;f_{z^{a+1} u^{c-1}}) 
	\\
&+ (a+1)(a-1)L'^{w'}_{z'}(z;f_{z^{a}\bar z u^{c-1}})
	+(a+1)aL'^{w'}_{z'}(z;f_{z^{a-1} u^{c}})\2L^w_{z\bar z} + \ldots,
	\\
g_{z\bar z^{a+1} u^{c-1}} + (& a+ 1)g_{\bar z^{a} u^{c}} \2L^w_{z\bar z}
	+ \2L^w_{z\bar z^{a+1} u^{c-1}} 
	\\
 = & \2L'^{w'}_{\bar z'z'^{a+1} u'^{c-1}} +
	L'^{w'}_{\bar z'}(\bar f_{\bar z^{a+1} u^{c-1}};z) + \ldots,
\end{split}
\end{equation}
where as before the dots stand for the terms already determined.
Furthermore, proceeding by induction on $a$,
we may also assume that the derivatives
$f_{z^{a-1} u^c}$, and therefore also $f_{z^{a}\bar z u^{c-1}}$
in view of \eqref{basic21},
are already determined
and will include them in the dot terms.

The derivatives $g_{z^{a} u^c}$, $a\ge2$, are
uniquely determined by the conditions
\begin{equation*}
\2L'^{w'}_{z'^a u'^c} =0,
\end{equation*}
or, equivalently,
\begin{equation*}
\2L'^{w'}_{z'^a u'^c} + \1{\2L'^{w'}_{\bar z'^a u'^c}} =0.
\end{equation*}
We shall now include $g_{z^{a} u^c}$ in the dot terms
assuming them being determined by induction on $a$. 
In case $c\ge 1$ we consider the sum and the difference
of the second identity in \eqref{bb0} divided by $(a+1)$ 
and the conjugate of the third identity there:
\begin{equation}\Label{}
\begin{split}
2g_{z^{a+1}\bar z u^{c-1}} &+ 
	\Big(\frac1{a+1}{\2L^w_{z^{a+1}\bar z u^{c-1}}} +
		\1{\2L^w_{z\bar z^{a+1} u^{c-1}}}\Big) 
	\\
	=& \Big(\frac1{a+1}\2L'^{w'}_{z'^{a+1}\bar z' u'^{c-1}} 
		+ \1{\2L'^{w'}_{z'^{a+1} u'^{c-1}}}\Big) + \ldots
	\\
\Big(\frac1{a+1} & \2L^w_{z^{a+1}\bar z u^{c-1}}  - 
	\1{\2L^w_{z\bar z^{a+1} u^{c-1}}}\Big) 
	\\
	= & \Big(\frac1{a+1}\2L'^{w'}_{z'^{a+1}\bar z' u'^{c-1}} - 
	\1{\2L'^{w'}_{z'z'^{a+1} u'^{c-1}}}\Big)
	+2L'^{w'}_{\bar z'}(\bar z;f_{z^{a+1} u^{c-1}}).
\end{split}
\end{equation}
Then the derivatives $g_{z^{a+1}\bar z u^{c-1}}$
 and $f_{z^{a+1} u^{c-1}}$
are uniquely determined by the normalization
\begin{equation}\Label{}
\2L'^{w'}_{z'^{a+1}\bar z' u'^{c-1}} \pm 
\1{\2L'^{w'}_{z'\bar z'^{a+1} u'^{c-1}}}
=0.
\end{equation}
By now we have determined all pure derivatives of $f$ of weight $k-1$
and shall include them in the dots.

Finally, for $a\ge b+2$, $b\ge2$, we have
\begin{equation}\Label{}
\begin{split}
ag_{z^{a}\bar z^b u^c} + 
ab g_{z^{a-1}\bar z^{b-1} u^{c+1}} \2L^w_{z\bar z}
+ \2L^w_{z^a\bar z^b u^c}
&= \2L'^{w'}_{z'^a\bar z'^b u'^c}+\ldots,\\
bg_{z^{b}\bar z^{a} u^c} + 
ab g_{z^{b-1}\bar z^{a-1} u^{c+1}}\2L^w_{z\bar z}
+ \2L^w_{z^{b}\bar z^{a} u^c}
&= \2L'^{w'}_{z'^{b}\bar z'^{a} u'^c} +\ldots,
\end{split}
\end{equation}
from where $g_{z^{a}\bar z^b u^c}$ is uniquely determined
by the condition
\begin{equation}\Label{n22}
\2L'^{w'}_{z'^a\bar z'^b u'^c} + \1{\2L'^{w'}_{z'^{b}\bar z'^{a} u'^c}}
=0.
\end{equation}
Thus we have determined all derivatives of $g$ of weight $k$.

Summarizing, we have determined all derivatives of $g$ of weight not
greater than $k$, and of $f$ of weight not greater than $k-1$.
We finally use \eqref{basic21} to determine the non-pure derivatives
$f_{z'^{a} \bar z'^b u'^c}$ of weight $k$ via the normalization
\begin{equation}\Label{n23}
\2L'^{\bar z'}_{z'^a \bar z'^b u'^c} = 0,
\end{equation}
and as the last step, determine all derivatives of $L'$ of weight $k-1$ from \eqref{basic}.

\subsection{Normalization summarized}
Collecting and simplifying our normalization \eqref{n100}, \eqref{n101}, \eqref{n22}
and \eqref{n23} we obtain a normal form for 
strongly nondegenerate almost CR structures:

\bt\Label{full}
For every formal power series 
$L(z,\bar z,u)\colon \C^n_z\to\C^n_{\bar z}\times\C_w$
without constant terms
corresponding to a strongly nondegenerate almost CR structure and every $r\in\R$, 
there exist unique formal power series 
$f(z,\bar z,u)\in\C^n$ and $g(z,\bar z,u)\in\R$
without constant and linear terms 
such that $g_{u^2}(0)=r$ and 
the map $h=\id + (f,g)$ transforms $L$ into $L'$
satisfying the normalization
\begin{equation}\Label{n100'}
\begin{split}
&\2L'^{\bar z'}_{z'^a \bar z'^b u'^c} = 0,\\
&\2L'^{w'}_{z'^a\bar z'^b u'^c} + 
\1{\2L'^{w'}_{z'^{b}\bar z'^{a} u'^c}}=0,\quad
\2L'^{w'}_{z'\bar z' u'^{c+1}} =0, \quad
\2L'^{w'}_{z'^{a+2}\bar z' u'^{c}} =0, \\
&\tr(\2L'^{w'}_{z'^2\bar z'^2 u'^c}) =0 ,\quad
\tr^2 (\2L'^{w'}_{z'^3\bar z'^2 u'^{c}})=0, \quad
\tr^3(\2L'^{w'}_{z'^3\bar z'^3 u'^{c}}) =0,
\end{split}
\end{equation}
for all $a,b,c\ge0$.
Furthermore, each partial derivative of $L'$, $f$, $g$ at $0$ is given by an universal polynomial in $r$ and (finitely many) derivatives of $L$.
\et

\section{Extrinsic normal form for hypersurface type almost CR structures}
Following the idea of quasi CR embeddings in \S\ref{quasi}
we refine the extrinsic normalization of Proposition~\ref{quasi-prop}
similarly to the intrinsic approach of \S\ref{hyper-form}.
As in \S\ref{quasi} we assume $L$ to be given and look for a real function $\phi(z,\bar z,u)$
normalized as in \eqref{phi-vanish} and \eqref{noharm} and a transformation $h=\id+(f,g)$
sending $L$ into $L'$ such that 
\eqref{phi-L} is satisfied with $L$ replaced by $L'$ along the Euler vector field \eqref{euler}.
To simplify the notation we shall drop $'$ from the variables
when indicating the derivatives and components of $L'$,
e.g.\ we shall write $L'^w_{zu}$ instead of $L'^{w'}_{z'u'}$.

Thus we evaluate \eqref{phi-L} 
with $L$ replaced by $L'$
 along the 
Euler vector field $e$:
\begin{equation}\Label{phi-Lexact}
(\2L'^{\bar z}, \2L'^w)= \big(0,  i (\id -i\phi_u)^{-1}\phi_z e\big).
\end{equation}
We differentiate \eqref{phi-Lexact} at $0$,
this time writing explicitly all terms of the maximum weight:
\begin{equation}\Label{der-lphi'}
\2L'^w_{z^{a}\bar z^b u^c} = 
ia (\phi_{z^{a}\bar z^b u^c} 
+ib\phi_{z^{a-1}\bar z^{b-1} u^{c+1}} \phi_{z\bar z})
+ R'_{a,b,c}(\phi_*),
\end{equation}
where $R'_{a,b,c}(\phi_*)$ is a polynomial in the derivatives $\phi_*$
of $\phi$ at $0$ of weight less than $a+b+2c$.
As before we assume \eqref{noharm} for $\phi$.
We obtain $\phi_{z\bar z}=\frac1{i}\2L'^w_{z\bar z}=\frac1{i}\2L^w_{z\bar z}$ and subsequently
\begin{equation}\Label{phi-L'}
i\phi_{z^{a}\bar z^b u^c} = \frac1{a} \2L'^{w}_{z^{a}\bar z^b u^c} 
-ib \phi_{z^{a-1}\bar z^{b-1} u^{c+1}} \2L^w_{z\bar z}
- \frac1{a} S'_{a,b,c}(\2L'^w_*),
\end{equation}
directly following from \eqref{der-lphi'} by induction, where
$S'_{a,b,c}(\2L'^w_*)$ is a polynomial in the derivatives
$\2L'^w_*$ of $\2L'^w$ at $0$ of weight less than $a+b+2c$.

As in \S\ref{quasi} we look for transformations \eqref{fg}
normalizing $\2L'^{\bar z}$ and $\phi$, where $\phi$ is a priori complex-valued and its reality  will be imposed as part of the normalization. We also assume 
the conditions \eqref{norm1} for both $L$ and $L'$
and therefore also \eqref{vanish-low}.
In particular, we have 
$L'_{z}=L_{z}$, $L'_{\bar z}=L_{\bar z}$ and
  $\2L'_{z\bar z}=\2L_{z\bar z}$.
We follow the strategy of \S\ref{hyper-form}
and assume for an integer $k\ge3$, that
we have determined all derivatives of $\phi$ and $g$ of weight less than $k$ and all
derivatives of $f$ of weight less than $k-1$.
As before, the dots will stand for the terms already determined.
The terms involving only $L$ are fixed 
and will also be included in the dots.

We first claim that in view of \eqref{phi-Lexact},
the normalization conditions
\begin{equation}\Label{n101}
\begin{split}
&\2L'^{w}_{z^a\bar z^b u^c} =0,\quad \min(a,b)\le 1, \quad (a,b,c)\ne (1,1,0), \\
&\tr(\2L'^{w}_{z^2\bar z^2 u^c}) =0 ,\quad
\tr^2 (\2L'^{w}_{z^3\bar z^2 u^{c}})=0, \quad
\tr^3(\2L'^{w}_{z^3\bar z^3 u^{c}}) =0,
\end{split}
\end{equation}
(which is a part of \eqref{n100'}), are equivalent to
\begin{equation}\Label{n102}
\begin{split}
&\phi_{z^a\bar z^b u^c} =0,\quad \min(a,b)\le 1, \quad (a,b,c)\ne (1,1,0), \\
&\tr(\phi_{z^2\bar z^2 u^c}) =0 ,\quad
\tr^2 (\phi_{z^3\bar z^2 u^{c}})=0, \quad
\tr^3(\phi_{z^3\bar z^3 u^{c}}) =0.
\end{split}
\end{equation}
Indeed, we assume either of \eqref{n101} and \eqref{n102} and prove the other set of conditions by 
the induction on the weight. Differentiating \eqref{phi-Lexact} as above and using
\eqref{n102} for terms of lower weight by either the assumption or the induction assumption, we obtain
$R'_{a,b,c}(\phi_*)=0$ in \eqref{der-lphi'} whenever either $\min(a,b)\le1$
or $\max(a,b)\le3$ and $(a,b)\ne(3,3)$. 
Since we have $\phi_{z^{a-1}\bar z^{b-1} u^{c+1}}=0$
for $(a,b)$ in this range, \eqref{der-lphi'} reduces to
$\2L'^w_{z^{a}\bar z^b u^c} = 
ia \phi_{z^{a}\bar z^b u^c}$, implying the induction step for those $(a,b,c)$.
In particular, we have $\tr(\phi_{z^2\bar z^2 u^c}) =0$
implying $\tr^3(\phi_{z^2\bar z^2 u^c}\phi_{z\bar z})=0$
(i.e.\ by the uniqueness of the decomposition \eqref{decomp}).
Then differentiating \eqref{phi-Lexact} $3$ times in each of $z$ and $\bar z$
and taking $\tr^3$ of both sides,
we conclude $\tr^3 (\2L'^{w}_{z^3\bar z^3 u^c})=\tr^3(\phi_{z^3\bar z^3 u^c})$,
completing the proof of the induction step and thus proving the claim.

We now continue our strategy following the lines of \S\ref{hyper-form} and consider \eqref{age4} for $a\ge4$. Assuming $g_{z^{a-1}\bar z^{a-1} u^{c+1}}$ and $\phi_{z^{a-1}\bar z^{a-1} u^{c+1}}$ being determined by induction on $a$ and using \eqref{phi-L'}, we obtain
$
i\phi_{z^a\bar z^a u^c}=g_{z^a\bar z^a u^c}+\ldots
$.
Then it is clear that the reality condition $\Im \phi_{z^a\bar z^a u^c}=0$ uniquely determines $g_{z^a\bar z^a u^c}$ and $\phi_{z^a\bar z^a u^c}$.

Finally, the normalization of the remaining terms is straightforward following the strategy of \S\ref{hyper-form} and is given by the remaining reality conditions of the form
$\phi_{z^{a}\bar z^b u^c}-\1{\phi_{z^{b}\bar z^a u^c}}=0$.
Summarizing we obtain:

\bt\Label{full-extrinsic}
For every formal power series 
$L(z,\bar z,u)\colon \C^n_z\to\C^n_{\bar z}\times\C_w$
without constant terms
corresponding to a strongly nondegenerate almost CR structure and every $r\in\R$,
there exist unique formal power series
$f(z,\bar z,u)\in\C^n$, $g(z,\bar z,u)\in\R$ and $\phi(z',\bar z',u')\in\R$ without constant and linear terms
satisfying the Chern-Moser normalization
\begin{equation}\Label{n100''}
\begin{split}
&\phi_{z'^a u'^c}=0,\quad
\phi_{z'\bar z' u'^{c+1}} =0, \quad
\phi_{z'^{a+2}\bar z' u'^{c}} =0, \\
&\tr(\phi_{z'^2\bar z'^2 u'^c}) =0 ,\quad
\tr^2 (\phi_{z'^3\bar z'^2 u'^{c}})=0, \quad
\tr^3(\phi_{z'^3\bar z'^3 u'^{c}}) =0
\end{split}
\end{equation}
for all $a,b,c\ge0$, and
such that $g_{u^2}(0)=r$ and the almost CR structure given by $L$ admits 
a quasi CR embedding at $0$ as the hypersurface $M'\subset\C^n_{z'}\times\C_{w'}$ given by
$$\Im w' = \phi(z',\bar z', \Re w')$$
via the map $(z,u)\mapsto (z+f(z,\bar z,u), u+g(z,\bar z,u))$, 
where $(z',u')$ are regarded as intrinsic coordinates on $M'$.
Furthermore, each partial derivative of $f$, $g$, $\phi$ at $0$ is given by an universal polynomial in $r$ and (finitely many) derivatives of $L$.
\et

Note that in view of Remark~\ref{forms-different}
the normal form in Theorem~\ref{full-extrinsic} is different from that of Theorem~\ref{full}.

\section{Uniqueness of the normal forms}
We here provide the details to the construction of extended adapted frames
described in \S\ref{intro} as well as the proof of Theorem~\ref{uniqueness}.

\subsection{$H$-equivalence}\Label{jets}
Let $J^k_{p,p'}(M,M')$ denote the space of all $k$-jets
of smooth maps between real manifolds $M$ and $M'$
with source $p\in M$ and target $p'\in M'$.
In local coordinates
a $k$-jet in $J^k(M,M')$ is represented by 
a $k$th order polynomial map.
Write $J^k(M,M')$ for the union of $J^k_{p,p'}(M,M')$ 
for all $p$ and $p'$.
Denote by $\pi^{k}_{k-1}\colon J^{k}(M,M')\to J^{k-1}(M,M')$
the natural projection.
We shall make use of the well-known fact that for each 
$\L\in J^{k-1}(M,M')$, the fiber $(\pi^{k}_{k-1})^{-1}(\L)$
has a canonical structure of an affine space modeled on 
the space of all symmetric $k$-linear maps from $T_pM\times\cdots\times T_pM$ 
into $T_{p'}M'$,
where $p$ and $p'$ are respectively the source and the target of $\L$.
The latter fact immediately follows from the chain rule.
Indeed, in local coordinates, points of $(\pi^{k}_{k-1})^{-1}(\L)$
are represented by the $k$th derivatives $D^k_p f$ of smooth maps 
$f\colon M\to M'$ representing the jets.
Then coordinate changes $\Phi$ and $\Psi$ on $M$ and $M'$ respectively
lead to the change 
\begin{equation}\Label{change}
D^k_p f\mapsto D^k_p(\Psi\circ f\circ\Phi^{-1}) 
= \Phi_* \circ  D^k_p  f \circ (\Phi^{-1}_*,\ldots, \Phi^{-1}_*)+ P(j^{k-1}_pf),
\end{equation}
where $P(j^{k-1}_pf)$ is a polynomial in the $(k-1)$-jet variables
(that depends on $\Phi$ and $\Psi$).
Thus one obtains an affine transformation rule for $D^k_p f$
when $j^{k-1}_pf$ is fixed as claimed.

Let now $M$ be an almost CR-manifold and consider the second jet space 
$J^2_{0,p}(\R,M)$ with fixed source 
$0\in \R$ and target $p\in M$.
Obviously the first jet space $J^1_{0,p}(\R,M)$
can be identified with the tangent space $T_pM$.
Then for a tanget vector $\xi\in T_pM $,
the fiber $(\pi^2_1)^{-1}(\xi) \subset J^2_{(0,p)}(\R,M)$
has the canonical structure of an affine space modeled on 
the space of all symmetric bilinear maps from $\R\times\R$ into $T_pM$,
which is canonically isomorphic to 
$T_pM$.
Hence for $\L_1,\L_2\in (\pi^2_1)^{-1}(\xi)$, 
we can consider the difference $\L_1-\L_2\in T_pM$.
We say that $\L_1$ and $\L_2$ are {\em $H_pM$-equivalent}
(or simply {\em $H$-equivalent}) if $\L_1-\L_2\in H_pM$.
This obviously defines an equivalence relation that we call
{\em $H$-equivalence}.

\subsection{Extended adapted frames}
Recall from \S\ref{intro}
that an {\em extended adapted frame} on $M$ at $p$ consists of a $\C$-basis $v_1,\ldots,v_n$ of $H_pM$,
a vector $v_{n+1}\in T_pM\setminus H_pM$ 
satisfying \eqref{orthog}
and an $H_pM$-equivalence class $[\L]\subset J^2_{0,p}(\R,M)$
with $\pi^1_0\L=v_{n+1}$.
The {\em standard extended adapted frame} on $\C^n\times\R$ at $0$
(with any almost CR structure satisfying $H_0=\C^n\times\{0\}$)
consists of the unit vectors $\frac{\d}{\d x_1},\ldots,\frac{\d}{\d x_n},\frac{\d}{\d u}$ and the $H_0M_0$-equivalence class represented by the linear curve $\gamma(t)=(0,t)\in \C^n\times \R$.

\bt\Label{unique-precise}
Let $L(z,\bar z,u)\colon \C^n_z\to \C^n_{\bar z}\times \C_w$
be a formal power series map without constant terms defining 
a strongly nondegenerate hypersurface type almost CR structure on $\C^n\times \R$ (in the formal sense) with a given extended adapted frame at $0$.
Then the following hold.
\begin{enumerate}
\item[(i)] 
There exists an unique formal map 
$h(z,\bar z,u)\in \C^n\times \R$ without constant terms
transforming the given extended adapted frame into the standard one
and $L$ into the normal form given by Theorem~\ref{full}.
\item[(ii)]
There exists an unique formal map
$h(z,\bar z,u)\in \C^n\times \R$  without constant terms
transforming the given extended adapted frame into the standard one
and realizing a quasi CR embedding of the given almost CR structure as a hypersurface 
$\Im w'=\phi(z',\bar z',\Re w)$ as in Theorem~\ref{full-extrinsic}.
\end{enumerate}
Furthermore, each partial derivative of $h$ and $L'$ corresponding to the normal form in each of {\rm(i)} and {\rm(ii)}
is given by a (universal) real-analytic function in the partial derivatives of $L$ and
the components of the given extended adapted frame. 
\et

\bpf
Let $v_1,\ldots,v_n,v_{n+1},[\L]$ be the given extended adapted frame.
Applying a linear transformation, we may assume that 
$v_1,\ldots,v_n,v_{n+1}$ is the standard adapted frame. 
Then the statements in (i) and (ii) are given by Theorems~\ref{full} and ~\ref{full-extrinsic}
respectively, where the class $[\L]$ uniquely determines the parameter $r=g_{u^2}(0)$.
\epf

\subsection{Proof of Theorem~\ref{system}}
Fix a point $p_0\in M$ and an extended adapted frame $\6F_p$ at each point $p$
in a neighborhood $U$ of $p_0$ in $M$. If $M$ is smooth (resp.\ real-analytic), we may choose
$U$ and $\6F_p$ that depends smoothly (resp.\ real-analytically) on $p$. 
In view of Theorem~\ref{unique-precise}, for every $p\in U$, 
the almost CR structure of $M$ at the reference point $p$ can be mapped via a map $h_p$ 
formally and uniquely
into its normal form $L'_p$ given by Theorem~\ref{full} such that
$\6F_p$ is mapped into the standard extended adapted frame at $0$.
Furthermore, the last statement of Theorem~\ref{unique-precise} implies
that the coefficients of $h_p$ and $L'_p$ depend smoothly (resp.\ real-analytically)
on $p$. In order to show the desired conclusion, we consider two cases.

{\bf Case 1}. The given $2$-jet $\L_0\in G^2_{p_0}(M,M')$ does not send the chosen
extended adapted frame at $p_0$ into any extended adapted frame in $M'$.
In that case $\L_0$ cannot be a $2$-jet of a CR-diffeomorphism
and furthermore a neighborhood $\Omega$ of $\L$ in $G^2(M,M')$ can be chosen
such that no $2$-jet in $\Omega$ is a $2$-jet of a CR-diffeomorphism.
Then the property $j^2_pf\in \Omega$ is never satisfied for a CR-diffeomorphism $f$
and hence the conclusion \eqref{sys} holds trivially with any choice of the map $\Phi$.

{\bf Case 2}. The given $2$-jet $\L_0\in G^2_{p_0}(M,M')$ sends the chosen
extended adapted frame at $p_0$ into an extended adapted frame in $M'$.
We first construct the map $\Phi$ only for $2$-jets $\L\in G^2_{p,p'}(M,M')$
such that $p\in U$ and $\L$ sends $\6F_p$
into some extended adapted frame $\6F_{p'}$ on $M'$ at $p'$.
Note that the set $S$ of such $\L$ is a smooth (resp.\ real-analytic) 
submanifold of $G^2(M,M')$.
Applying again Theorem~\ref{unique-precise} for the almost CR structure
of $M'$ at the reference point $p'$ and the frame $\6F_{p'}$,
we obtain an uniquely determined map $h_{\L}$, depending on $\L$, sending $M'$ into its normal form $M'_{\L}$ such that $\6F_{p'}$ is sent to the standard extended adapted frame.
Then, for any $p\in U$ and CR-diffeomorphism $f$ between open pieces of $M$ and $M'$,
defined in a neighborhood of $p$, the formal maps $h_p$ and $h_{j^2_pf}\circ \2f$
must coincide by the uniqueness in Theorem~\ref{unique-precise}.
Here $\2f$ denotes the formal map given by the Taylor series of $f$ at $p$.
Hence $\2f=h_{j^2_pf}^{-1}\circ h_p$ and,
in particular, $j^3_pf =j^3_p\2f = j^3_p(h_{j^2_pf}^{-1}\circ h_p)$.
We then set $\Phi(\L):=j^3_p(h_{\L}^{-1}\circ h_p)$.
It remains to choose any open neighborhood $\Omega$ of $\L$ in $G^2(M,M')$
with a smooth (resp.\ real-analytic) retraction $r\colon \Omega\to S\cap\Omega$
and replace $\Phi$ by $\Phi\circ r$, now defined in $\Omega$.
The proof is complete.


\begin{thebibliography}{BER96bb}

\bibitem[BER99]{BERbook} {\bf Baouendi,~M.S.; Ebenfelt,~P.; Rothschild,~L.P.} --- {\em Real Submanifolds in Complex Space and Their Mappings}. Princeton Math. Series {\bf 47}, Princeton Univ. Press, 1999.

\bibitem[BRWZ04]{BRWZ} {\bf Baouendi,~M.S.; Rothschild,~L.P.; Winkelmann,~J.; Zaitsev,~D.} --- Lie group structures on groups of diffeomorphisms and applications to CR manifolds.  {\em Ann. Inst. Fourier (Grenoble)}  {\bf 54}  (2004),  no. 5, 1279--1303.

\bibitem[Ca32]{CaE} {\bf Cartan,~E.} --- Sur la g\'eom\'etrie pseudo-conforme des hypersurfaces de deux variables complexes, I. {\em Ann. Mat. Pura Appl.} {\bf 11} (1932), 17--90. ({\em \OE uvres compl\`etes}, Part. II, Vol. 2, Gauthier-Villars, 1952, 1231--1304); II. {\em Ann. Sc. Norm. Sup. Pisa} {\bf 1} (1932), 333--354. ({\em \OE uvres compl\`tes}, Part. III, Vol. 2, Gauthier-Villars, 1952, 1217--1238).


\bibitem[CM74]{CM} {\bf Chern,~S.S; Moser,~J.K.} --- Real hypersurfaces in complex manifolds. {\em Acta Math.} {\bf 133} (1974), 219--271.

\bibitem[P07]{Po} {\bf Poincar\'e,~H.} --- Les fonctions analytiques de deux variables et la repr\'esentation conforme. {\em Rend.~Circ.~Mat.~Palermo} {\bf 23}, 185--220, (1907).


\bibitem[E01]{Ejets} {\bf Ebenfelt,~P.} --- Finite jet
determination of holomorphic mappings at the
boundary. {\em Asian J. Math.} {\bf 5} (2001),
637-662.

\bibitem[KZ05]{KZ} {\bf Kim,~S.-Y.; Zaitsev,~D.} --- 
Equivalence and embedding problems for
CR-structures of any codimension. {\em Topology}, {\bf 44} (3), (2005), 557--584.

\bibitem[Kr98a]{Kr1} {\bf Kruglikov,~B.S.} ---
Nijenhuis tensors and obstructions to the construction of pseudoholomorphic mappings. 
{\em Mat. Zametki} {\bf 63} (1998), no. 4, 541--561; translation in
{\em Math. Notes} {\bf 63} (1998), no. 3-4, 476--493.

\bibitem[Kr98b]{Kr2} {\bf Kruglikov,~B.S.} --- On some classification problems in four-dimensional geometry: distributions, almost complex structures, and the generalized Monge-Ampre equations. 
{\em Mat. Sb.} {\bf 189} (1998), no. 11, 61--74; translation in
{\em Sb. Math.} {\bf 189} (1998), no. 11-12, 1643--1656.

\bibitem[Ta62]{Ta62} {\bf Tanaka,~N.} --- On the pseudo-conformal geometry of hupersurfaces of the space of $n$ complex variables. {\em J. Math.\ Soc.\ Japan} {\bf 14} (1962), 397--429.

\bibitem[To08]{To} {\bf Tonejc, J.} --- Normal forms for almost complex structures.
{\em Internat. J. Math.} {\bf 19} (2008), no. 3, 303--321. 

\end{thebibliography}
\end{document}